\documentclass{amsart}
\usepackage{enumerate}
\usepackage{amssymb}

\usepackage{color}

\newtheorem{theorem}{Theorem}
\newtheorem{corollary}[theorem]{Corollary}
\newtheorem{definition}[theorem]{Definition}
\newtheorem{example}[theorem]{Example}
\newtheorem{lemma}[theorem]{Lemma}

\newtheorem{remark}[theorem]{Remark}

\begin{document}

\title{Commutator estimates in $W^*$-algebras}
\author{A. F. Ber}
\address{Department of Mathematics, Tashkent State University, Uzbekistan }
\email{ber@ucd.uz}

\author{F. A. Sukochev}
\address{School of Mathematics and Statistics, University of New South Wales, Sydney, NSW 2052, Australia }
\email{f.sukochev@unsw.edu.au}

\maketitle



\bigskip
\begin{abstract}
\noindent Let $\mathcal{M}$ be a $W^*$-algebra and let
$LS(\mathcal{M})$ be the algebra of all locally measurable operators
affiliated with $\mathcal{M}$. It is shown that for any
self-adjoint element $a\in LS(\mathcal{M})$ there exists a
self-adjoint element $c_{_{0}}$ from the center of
$LS(\mathcal{M})$, such that for any $\varepsilon>0$ there exists
a unitary element $ u_\varepsilon$ from $\mathcal{M}$, satisfying
$|[a,u_\varepsilon]| \geq (1-\varepsilon)|a-c_{_{0}}|$. A
corollary of this result is that for any derivation $\delta$ on
$\mathcal{M}$ with the range in a (not necessarily norm-closed)
ideal $I\subseteq\mathcal{M}$, the derivation $\delta$ is inner,
that is $\delta(\cdot)=\delta_a(\cdot)=[a,\cdot]$, and $a\in I$.
Similar results are also obtained for inner derivations on
$LS(\mathcal{M})$.
\end{abstract}



\section{Introduction}

Let $\mathcal{M}$ be a $W^*$-algebra and let $Z(\mathcal{M})$ be
the center of $\mathcal{M}$. Fix $a\in\mathcal{M}$ and consider
the inner derivation $\delta_a$ on $\mathcal{M}$ generated by the
element $a$, that is $\delta_a(\cdot):=[a,\cdot]$. Obviously,
$\delta_a$ is a linear bounded operator on
$(\mathcal{M},\|.\|_{\mathcal{M}})$, where $\|.\|_{\mathcal{M}}$
is a $C^*$-norm on $\mathcal{M}$. It is well known (see
e.g. \cite[Theorem 4.1.6]{Sak}) that there exists $c\in
Z(\mathcal{M})$ such that the following estimate holds:
$\|\delta_a\|\geq\|a-c\|_{\mathcal{M}}$. In view of this result,
it is natural to ask whether there exists an element
$y\in\mathcal{M}$ with $\|y\|\leq 1$ and $c\in Z(\mathcal{M})$
such that $|[a,y]|\geq|a-c|$?



The following estimate easily follows from the main result of the present article:
for every  self-adjoint element $a\in\mathcal{M}$ there exists an
element $c\in Z(\mathcal{M})$ and the family $\{u_\varepsilon\}_{\varepsilon>0}$ of unitary operators from $\mathcal{M}$ such that
\begin{equation}\label{zero main}
|\delta_a(u_\varepsilon)|\geq
(1-\varepsilon)|a-c|,\ \forall \varepsilon>0.
\end{equation}

The estimate above is actually sharp and, with its aid, we shall
easily show (see Corollary \ref{c1} below) that every derivation
$\delta$ on $\mathcal{M}$ taking its values in a (not necessary
closed in the norm $\|.\|_{\mathcal{M}}$) two-sided ideal
$I\subset\mathcal{M}$ has the form $\delta=\delta_a$, where $a\in
I$. This result can be further reformulated in two equivalent
forms (Corollary \ref{c1_1} and \ref{c1_2}) yielding
generalizations and complements to classical results of J. Calkin
\cite{Cal} and M.J. Hoffman \cite{Ho} obtained originally for the
special case when $\mathcal{M}$ coincides with the algebra $B(H)$
of all bounded linear operators on a Hilbert space $H$.
It should be pointed out that our approach to the proof of
Corollaries \ref{c1_1} and \ref{c1_2} is based on the estimate
\eqref{zero main} and appears to be more direct than those
employed in \cite{Cal} and \cite{Ho}. In Section
\ref{applications} below we present a number of other extensions
of results from \cite{Cal, Ho}, in particular to ideals of
measurable (unbounded) operators, which significantly extend
recent results from \cite{AAK, BdPS} obtained under an additional
assumption that $\mathcal{M}$ is of type $I$.

Further, we recall the following well-known problem, which is also somewhat relevant to our discussion of
derivations taking values in ideals.

Let $\mathcal{N}$ be a von Neumann subalgebra of the von Neumann
algebra $\mathcal{M}$ and  let $I$ be an arbitrary (two-sided) ideal in
$\mathcal{M}$. What conditions should  be imposed on
$\mathcal{M},\ \mathcal{N},\ I$ to guarantee that for every
derivation $\delta: \mathcal{N}\rightarrow I$ there exists $a\in
I$ so that the equality $\delta=\delta_a$ holds?

Different partial solutions of this problem can be found
\cite{JP,KW,PR}. In the present paper we present a positive
solution in the special case when $\mathcal{M}$ is an arbitrary
von Neumann algebra and when $\mathcal{N}=\mathcal{M}$. Let us
note that if, in addition, we assume that the ideal $I$ is closed
with respect to the norm $\|.\|_{\mathcal{M}}$, then the positive
solution can be obtained directly from the Dixmier approximation
theorem (see e.g. \cite[Theorem 2.1.16]{Sak}). However, the latter
theorem is inapplicable when the ideal $I$ is not closed in the
$C^*$-norm on $\mathcal{M}$.

Analogous (but much harder) questions to those discussed above can
be also reformulated in a more general setting of the theory of
non-commutative integration initiated by I.E. Segal \cite{Se} (for
alternative approach to this theory, see E. Nelson's paper
\cite{Ne}). In these reformulations (see e.g. \cite{BdPS, BS1}),
the  $W^*$-algebra $\mathcal{M}$ is replaced with a larger algebra
of \lq measurable' operators affiliated with $\mathcal{M}$ and the
ideal $I$ in $\mathcal{M}$ is replaced with an ideal of measurable
operators. The most general algebra considered in the theory of
noncommutative integration to date is the classical algebra
$LS(\mathcal{M})$ (see \cite{San}) of all locally measurable
operators affiliated with $\mathcal{M}$ (all necessary definitions
are given in Section \ref{preliminaries} below). It is important
to emphasize that our methods are totally different from the
methods employed in \cite{Cal,Ho,JP,KW,PR,AAK,BdPS} and are strong enough
to enable us (see Theorem \ref{main} and Corollaries
\ref{c2},\ref{c3} below) to resolve all these questions also in
the setting of the algebra $LS(\mathcal{M})$. A number of such
applications to symmetric spaces of measurable operators (some of
them are in spirit of papers \cite{Cal,Ho,KW,AAK,BdPS}) are presented in
Section \ref{applications}.

Our main result in this paper is the following theorem.

\begin{theorem}
\label{main}
Let $\mathcal{M}$ be a $W^*$-algebra and let $a=a^*\in LS(\mathcal{M})$.
\begin{enumerate}[(i)]
\item If $\mathcal{M}$ is a finite $W^*$-algebra or else a purely
infinite $\sigma$-finite $W^*$-algebra, then there exists
$c_{_{0}}=c_{_{0}}^*\in Z(LS(\mathcal{M}))$ and
$u_{_{0}}=u_{_{0}}^*\in U(\mathcal{M})$, such that
\begin{equation}\label{first main}
|[a,u_{_{0}}]|=u_{_{0}}^*|a-c_{_{0}}|u_{_{0}} + |a-c_{_{0}}|,
\end{equation}
where $U(\mathcal{M})$ is the group of all unitary elements in
$\mathcal{M}$; \item There exists $c_{_{0}}=c_{_{0}}^*\in
Z(LS(\mathcal{M}))$, so that for any $\varepsilon > 0$ there
exists $u_\varepsilon=u_\varepsilon^* \in U(\mathcal{M})$ such
that
\begin{equation}\label{second main}
|[a,u_\varepsilon]| \geq (1-\varepsilon)|a-c_{_{0}}|.
\end{equation}
\end{enumerate}
\end{theorem}

The main ideas of the proof of Theorem \ref{main} are firstly
demonstrated in Section \ref{the outline} and then presented in
full in Section \ref{the proof}. The proof of Theorem \ref{main}
in the special case when $\mathcal{M}$ is a factor can be also
found in our earlier paper \cite{BS1}. We point out that the proof
in the general case has required a significant new technical
insight and that our direct approach may be of interest in its own
right (not in the least due to the fact that we do not use any
complicated technical devises such as Boolean valued analysis and
direct integrals).

\begin{remark}\label{trivial}Observe that the equality \eqref{first main} trivially yields
the estimate \eqref{second main} even for the case
$\varepsilon=0$. Nevertheless, the result of Theorem
\ref{main}(ii) is still sharp. Indeed, if $\mathcal{M}$ is an
infinite semifinite $\sigma$-finite factor, then there exists a
self-adjoint element $a\in LS(\mathcal{M})$ such that for every
$\lambda\in\mathbb{C}$ and $u\in U(\mathcal{M})$ the inequality
$|[a,u]|\geq |a-\lambda\mathbf{1}|$ fails \cite{BS1}. Hence, the
multiplier $(1-\varepsilon)$ in the part $(ii)$ of Theorem
\ref{main} can not be omitted.
\end{remark}

\section{Preliminaries}\label{preliminaries}

For details on von Neumann algebra theory, the reader is referred
to e.g. \cite{Dix}, \cite{KR2}, \cite{Sak} or \cite{Tak}. General
facts concerning measurable operators may be found in \cite{Ne},
\cite{Se} (see also \cite[Chapter IX]{Ta2}). For the convenience
of the reader, some of the basic definitions are recalled.

Let $\mathcal{M}$ be a von Neumann algebra on a Hilbert space $H$.
The set of all self-adjoint projections (respectively, all unitary
elements) in $\mathcal{M}$ is denoted by
$P\left(\mathcal{M}\right)$ (respectively,
$U\left(\mathcal{M}\right)$). We use the notation
$\mathbf{s}(x),\mathbf{l}(x),\mathbf{r}(x),\mathbf{c}(x)$ to
denote the support, left support, right support, central support
respectively of an element $x\in\mathcal{M}$.

Let $p,q\in P\left(\mathcal{M}\right)$. The projections
$p$ and $q$ are said to be equivalent, if there exists
a partial isometry $v\in\mathcal{M}$, such that $v^*v=p,\ vv^*=q$. In this case, we write $p\sim q$.
The fact that the projections $p$ and $q$ are not equivalent is recorded as $p\nsim q$.
If there exists a projection $q_1\in P\left(\mathcal{M}\right)$ such that $q_1\leq p,\ q_1\sim q$, then we write
$q\preceq p$. If $q\preceq p$ and $p\nsim q$, then we employ the notation $q\prec p$.

A linear operator $x:\mathfrak{D}\left( x\right) \rightarrow
H $, where the domain $\mathfrak{D}\left( x\right) $ of $x$ is a linear
subspace of $H$, is said to be {\it affiliated} with $\mathcal{M}$ if $yx\subseteq
xy$ for all $y\in \mathcal{M}^{\prime }$ (which is denoted by $x\eta
\mathcal{M}$). A linear operator $x:\mathfrak{D}\left( x\right) \rightarrow
H $ is termed {\it measurable} with respect to $\mathcal{M}$ if $x$ is closed,
densely defined, affiliated with $\mathcal{M}$ and there exists a sequence $%
\left\{ p_{n}\right\} _{n=1}^{\infty }$ in $P\left( \mathcal{M}\right) $
such that $p_{n}\uparrow \mathbf{1}$, $p_{n}\left( H\right) \subseteq
\mathfrak{D}\left( x\right) $ and $p_{n}^{\bot }$ is a finite projection
(with respect to $\mathcal{M}$) for all $n$. It should be noted that the
condition $p_{n}\left( H\right) \subseteq \mathfrak{D}\left( x\right) $
implies that $xp_{n}\in \mathcal{M}$. The collection of all measurable
operators with respect to $\mathcal{M}$ is denoted by $S\left( \mathcal{M}%
\right) $, which is a unital $\ast $-algebra with respect to strong sums and
products (denoted simply by $x+y$ and $xy$ for all $x,y\in S\left( \mathcal{M%
}\right) $).

Let $a$ be a self-adjoint operator affiliated with $\mathcal{M}$.
We denote its spectral measure by $\{e^a\}$. It is known if $x$ is
a closed operator in $H$ with the polar decomposition $x = u|x|$
and $x\eta \mathcal{M}$, then $u\in\mathcal{M}$ and $e\in
\mathcal{M}$ for all projections $e\in \{e^{|x|}\}$. Moreover,
$x\in S(\mathcal{M})$ if and only if $x$ is closed, densely
defined, affiliated with $\mathcal{M}$ and $e^{|x|}(\lambda,
\infty)$ is a finite projection for some $\lambda> 0$. It follows
immediately that in the case when $\mathcal{M}$ is a von Neumann
algebra of type $III$ or a type $I$ factor, we have
$S(\mathcal{M})= \mathcal{M}$. For type $II$ von Neumann algebras,
this is no longer true.

A operator $x\in S\left( \mathcal{M}\right) $ is called $\tau
$-{\it measurable} if there exists a sequence
$\left\{p_{n}\right\} _{n=1}^{\infty }$ in $P\left(
\mathcal{M}\right) $ such that $p_{n}\uparrow \mathbf{1}$,
$p_{n}\left( H\right) \subseteq \mathfrak{D}\left( x\right) $ and
$\tau \left( p_{n}^{\bot }\right) <\infty $ for all $n$. The
collection $S\left( \tau \right) $ of all $\tau $-measurable
operators is a unital $\ast $-subalgebra of $S\left(
\mathcal{M}\right) $ denoted by $S\left( \mathcal{M}, \tau\right)
$. It is well known that a linear operator $x$ belongs to $S\left(
\mathcal{M}, \tau\right) $ if and only if $x\in S(\mathcal{M})$
and there exists $\lambda>0$ such that $\tau(e^{|x|}(\lambda,
\infty))<\infty$.

A closed operator $x\eta\mathcal{M}$ is called locally
measurable if there exists a sequence $\left\{z_{n}\right\}
_{n=1}^{\infty }$ of central projections in $\mathcal{M}$ such
that $z_{n}\uparrow \mathbf{1}$ and $xz_{n}\in S(\mathcal{M})$ for
any $n\in \mathbb{N}$. The collection of all locally measurable
operators with respect to $\mathcal{M}$ is denoted by $LS\left(
\mathcal{M} \right)$, which is a unital $\ast$-algebra with
respect to strong sums and products (denoted simply by $x+y$ and
$xy$ for all $x,y\in LS\left( \mathcal{M}\right) $).

It follows directly from the definition of local measurability,
that every measurable operator $c$ with respect to
$Z(\mathcal{M})$ is a locally  measurable operator and so
$S(Z(\mathcal{M}))\subset Z(LS(\mathcal{M}))$. On the other hand,
if $c\in Z(LS(\mathcal{M}))$ and $c=v|c|$ is a polar decomposition
for  $c$, then $v\in Z(\mathcal{M})$ and the spectral family
$e^{|c|}_\lambda$ of $|c|$ belongs to $Z(\mathcal{M})$ as well.
This means that $|c|\in S(Z(\mathcal{M}))$. Hence, $c=v|c|\in
S(Z(\mathcal{M}))$. So, we have
$$Z(LS(\mathcal{M}))=S(Z(\mathcal{M})).$$
In particular, $*$-algebra $Z(LS(\mathcal{M}))$ is $*$-isomorphic
to the commutative $*$-algebra $L^0(\Omega,\Sigma,\mu)$ of all
measurable complex valued functions on a measurable space with a
locally finite measure \cite{Se}.

Suppose now that two von Neumann algebras $\mathcal{M}$ and
$\mathcal{N}$ are $*$-isomorphic. In this case, it follows from
\cite{Ch} that the algebras $LS(\mathcal{M})$ and
$LS(\mathcal{N})$ of locally measurable operators affiliated with
von Neumann algebras $\mathcal{M}$ and $\mathcal{N}$ respectively
are also $*$-isomorphic. Assuming now that $\mathcal{M}$ is an
arbitrary $W^*$-algebra and $\mathcal{N}$ is a von Neumann algebra
$*$-isomorphic with $\mathcal{M}$, we see that the algebra
$LS(\mathcal{N})$ is uniquely defined up to $*$-isomorphism. In
view of this observation, the algebra $LS(\mathcal{M})$ of all
locally measurable operators affiliated with an arbitrary
$W^*$-algebra $\mathcal{M}$ is well defined.

\section{Applications to
derivations and essential commutants}\label{applications}

Recall that a \textit{derivation} on a complex algebra $A$ is a linear map
$\delta :A\rightarrow A$ such that
\begin{equation*}
\delta \left( xy\right) =\delta \left( x\right) y+x\delta \left(
y\right) ,\ \ \ x,y\in A.
\end{equation*}
If $a\in A$, then the map $\delta _{a}:A\rightarrow A$, given by
$\delta _{a}\left( x\right) =\left[ a,x\right] $, $x\in A$, is a
derivation. A derivation of this form is called \textit{inner}. In
this section we demonstrate a number of corollaries from Theorem
\ref{main} to inner derivations on the algebra $A$ (in the setting
when $A=\mathcal{M}$ or $A=LS(\mathcal{M})$) taking value in some
two-sided ideal in $A$.

\subsection{ Applications to ideals in $W^*$-algebras}

{Our results here extend those in \cite{Cal, Ho}. We are also
motivated by results in \cite{JP, KW, PR}.

\begin{corollary}
\label{c1} Let $\mathcal{M}$ be a $W^*$-algebra and let $I$ be an
ideal in $\mathcal{M}$. Let $\delta: \mathcal{M} \rightarrow I$ be
a derivation. Then there exists an element $a\in I$, such that
$\delta=\delta_a=[a,\cdot]$.
\end{corollary}
\begin{proof}
Since $\delta$ is a derivation on a $W^*$-algebra, it is necessarily inner
\cite[Theorem 4.1.6]{Sak}. Thus, there exists an element $d\in\mathcal{M}$, such that
$\delta(\cdot)=\delta_d(\cdot)=[d,\cdot]$. It follows from our hypothesis that
$[d,\mathcal{M}]\subseteq I$.

Using \cite[Lemma 7]{BS1}, we obtain
$[d^*,\mathcal{M}]=-[d,\mathcal{M}]^*\subseteq I^*=I$ and
$[d_k,\mathcal{M}]\subseteq I,\ k=1,2$, where $d=d_1+i d_2$,
$d_k=d_k^*\in \mathcal{M}$, for $k=1,2$. It follows now from
Theorem \ref{main},  that there exist $c_1,c_2\in Z(\mathcal{M})$
and $u_1,u_2\in U(\mathcal{M})$, such that $|[d_k,u_k]|\geq
1/2|d_k-c_k|$ for $k=1,2$. Again applying \cite[Lemma 7]{BS1}, we
obtain $d_k-c_k\in I$, for $k=1,2$. Setting
$a:=(d_1-c_1)+i(d_2-c_2)$, we deduce that $a\in I$ and
$\delta=[a,\cdot]$.
\end{proof}

Corollary \ref{c1} can be restated as follows.

\begin{corollary}
\label{c1_1} Let $\mathcal{M}$ be a $W^*$-algebra, let $I$ be an
ideal in $\mathcal{M}$ and let $\pi:
\mathcal{M}\rightarrow\mathcal{M}/I$ be a canonical epimorphism.
Then $\pi^{-1}(center(\mathcal{M}/I))=Z(\mathcal{M})+I$
\end{corollary}
\begin{proof}
Let $a\in \pi^{-1}(center(\mathcal{M}/I))$. Then $[a,x]=ax-xa\in I$ for any $x\in\mathcal{M}$.
Therefore the inner derivation $\delta_a$ satisfies the condition of Corollary \ref{c1}.
Then there exists an element $c\in Z(\mathcal{M})$ such that $a+c\in I$.
Therefore $a\in Z(\mathcal{M})+I$, that is
$\pi^{-1}(center(\mathcal{M}/I))\subset Z(\mathcal{M})+I$. The converse inclusion is trivial.
\end{proof}

In the special case when $\mathcal{M}=B(H)$, where $H$ is a
separable Hilbert space the result of Corollary \ref{c1_1}
coincides with that of \cite[Theorem 2.9]{Cal}. Thus, we can view
Corollary \ref{c1_1} as an extension of the classical Calkin
theorem to arbitrary $W^*$-algebras. This result has the following
~\lq\lq lifting\rq\rq~ interpretation. Let $\phi$ be an
epimorphism from $W^*$-algebra $\mathcal{M}$ on an arbitrary
algebra $\mathcal{A}$ and let $a$ be an element from the center of
the algebra $\mathcal{A}$. Then there exist an element $z\in
Z(\mathcal{M})$ such that $\phi(z)=a$.

Let us give another important corollary of Theorem \ref{main}
related to the notions of essential commutant and multiplier
ideals \cite{Ho}. To this end, let us fix a $W^*$-algebra
$\mathcal{M}$ and two self-adjoint ideals $I,J$ in $\mathcal{M}$. We set
$$
I:J=\{x\in\mathcal{M}:\ xJ\subset I\}
$$
and
$$
D(J,I)=\{x\in\mathcal{M}:\ [x,y]\in I,\ \forall y\in J\}.
$$
Observe that  $I:J$ is an ideal in $\mathcal{M}$. In particular,
$(I:J)^*=I:J=\{x\in\mathcal{M}:\ Jx\subset I\}$.

\begin{corollary}\label{c1_2} For any $W^*$-algebra
$\mathcal{M}$ and any ideals $I,J$ in $\mathcal{M}$ we have
$$D(J,I) = I:J+Z(\mathcal{M}).$$
\end{corollary}
\begin{proof}

Let $a\in I:J$, $c\in Z(\mathcal{M})$, $b\in J$.
Then $[a+c,b]=[a,b]=ab-ba\in I-I\subset I$, that is $a+c\in D(J,I)$.
Therefore $I:J+Z(\mathcal{M})\subset D(J,I)$.

In order to prove that $D(J,I)\subset I:J+Z(\mathcal{M})$, fix an
element $a\in D(J,I)$. For an arbitrary $x\in J, y\in
\mathcal{M}$, we have $x[a,y]+[a,x]y=xay-xya+axy-xay=[a,xy]\in
[a,J] \subset I,\ [a,x]y\in Iy\subset I$. Hence, $x[a,y]\in I$.
Since $x$ is an arbitrary element from $J$, we obtain $[a,y]\in
I:J$ and so $[a,\mathcal{M}]\subset I:J$. Now, it follows from
Corollary \ref{c1} that $a+c\in I:J$ for some $c\in
Z(\mathcal{M})$. Consequently, $a=(a+c)-c\in I:J+Z(\mathcal{M})$.
\end{proof}

In the special case, when $\mathcal{M}=B(H)$, the result of
Corollary \ref{c1_1} coincides with that of \cite[Theorem
1.1]{Ho}. The proof given there does not extend to arbitrary von
Neumann algebras and is spelt out only for the case of a separable
Hilbert space $H$ (the case of a nonseparable $H$ requires a
substantial effort outlined in \cite{Ho}). The proof given above
works for an arbitrary $W^*$-algebra $\mathcal{M}$ and hence for
an arbitrary von Neumann algebra represented on an arbitrary
Hilbert space. 

Classical examples of normed ideals $I$ satisfying the assumptions
of Corollary \ref{c1} above are given by symmetric operator ideals
\cite{GK1, GK2, Schatten, Simon}.

\begin{definition}\label{opideal}
A linear subspace ~$\mathcal I$~ in the von Neumann algebra
$\mathcal M$ equipped with a norm $\|\cdot\|_{\mathcal I}$ is said
to be {\it a symmetric operator ideal} if\\
(1) $\| S\|_{\mathcal I}\geq \| S\|$ for all $S\in \mathcal I$,\\
(2) $\| S^*\|_{\mathcal I} = \| S\|_{\mathcal I}$ for all $S\in \mathcal I$,\\
(3) $\| ASB\|_{\mathcal I}\leq \| A\| \:\| S\|_{\mathcal I}\| B\|$ for all
$S\in \mathcal I$, $A,B\in \mathcal M$.\\
\end{definition}

Observe, that every symmetric operator ideal ~$\mathcal I$~ is a
two-sided ideal in $\mathcal M$, and therefore by \cite[I.1.6,
Proposition 10]{Dix}, it follows from $0\leq S \leq T$ and $T \in
{\mathcal I}$ that $S \in {\mathcal I}$ and $\| S\|_{\mathcal I}
\leq \| T\|_{\mathcal I}$.

\begin{corollary}
\label{c11} Let $\mathcal{M}$ be a $W^*$-algebra, let $I$ be a
symmetric operator ideal in $\mathcal{M}$ and let $\delta:
\mathcal{M} \rightarrow I$ be a self-adjoint derivation. Then
there exists an element $a\in I$, satisfying the inequality
$\|a\|_I\leq\|\delta\|_{\mathcal{M}\rightarrow I}$ and such that
$\delta=\delta_a=[a,\cdot]$.
\end{corollary}
\begin{proof} Firstly, we observe that
$\|\delta\|_{\mathcal{M}\rightarrow I}<\infty$. Indeed, we have
$\delta=\delta_a,\ a\in I$ and therefore $\|\delta(x)\|_I=\|ax-xa\|_I\leq
\|ax\|_I+\|xa\|_I\leq 2\|a\|_I \|x\|_{\mathcal{M}}$, that is
$\|\delta\|_{\mathcal{M}\rightarrow I}\leq 2\|a\|_I<\infty$.

Let now $\delta$ be a self-adjoint derivation on $\mathcal{M}$,
that is $\delta(\cdot)=\delta_d(\cdot)=[d,\cdot]$ for some $d\in
\mathcal{M}$, such that $[d,x]^*=[d,x^*]$ for all
$x\in\mathcal{M}$. We have $x^*d^*-d^*x^*=dx^*-x^*d$, that is,
$x^*(d^*+d)=(d^*+d)x^*$ for all $x\in\mathcal{M}$. This
immediately implies $\mathfrak{Re}(d)\in Z_h(\mathcal{M})$ and so,
we can safely assume that $\delta=\delta_{i\cdot d}=[i\cdot
d,\cdot ]$, where $d$ is a self-adjoint operator from
$\mathcal{M}$. Fix $\varepsilon>0$ and let $c_0\in
Z_h(\mathcal{M})$, $u_\varepsilon \in U(\mathcal{M})$ be such that
$$|[d,u_\varepsilon]| \geq (1-\varepsilon)|d-c_0|.
$$ The assumption on $(I,\|\cdot\|)$ guarantees that
$(1-\varepsilon)\|d-c_0\|_I\leq
\|\delta(u_\varepsilon)\|_I\leq \|\delta\|_{\mathcal{M}\rightarrow
I}$. Since $\varepsilon$ was chosen arbitrarily, we conclude that
$\|d-c_0\|_I\leq
\|\delta\|_{\mathcal{M}\rightarrow I}$. Setting
$a=i(d-c_0)$ completes the proof.
\end{proof}

If the von Neumann algebra $\mathcal M$ is equipped with a
faithful normal semi-finite trace $\tau$, then the set
$$\mathcal{L}_p(\mathcal{M})=\left\{S\in {\mathcal{M}}:\ \tau(|S|^p)<\infty\right\}$$
equipped  with a standard norm
$$\|S\|_{\mathcal{L}_p(\mathcal{M})}=\max\{\|S\|_{B(H)},\ \tau(|S|^p)^{1/p}\}$$
is called Schatten-von Neumann $p$-class. In the type $I$ setting
these are the usual Schatten-von Neumann ideals of compact
operators \cite{GK1, GK2, Schatten, Simon}. The result of
Corollary \ref{c11} complements results given in \cite[Section
6]{KW}, where derivations from some subalgebras of $\mathcal{M}$
into Schatten-von Neumann $p$-classes were studied.


\subsection{Applications to ideals in $LS(\mathcal{M})$}

We begin by proving an analogue of Corollary \ref{c1} for ideals
of (unbounded) locally measurable operators. The following result
significantly strengthens  \cite[Proposition 6.17]{BdPS} where a
similar assertion was established under additional assumptions
that $\mathcal{M}=L_\infty (\nu)\bar \otimes B(H)$ (here,
$L_\infty (\nu)$ is an algebra of all bounded measurable functions
on a measure space) and $\mathcal{A}$ is an absolutely solid
algebra such that $\mathcal{M}\subset\mathcal{A}$. Similarly, the
result below complements the main result of \cite{AAK}, which
considered the case of an arbitrary von Neumann algebra
$\mathcal{M}$ of type $I$ and algebras
$\mathcal{A}=S(\mathcal{M},\tau),S(\mathcal{M}), LS(\mathcal{M})$.
Our approach here is completely different from techniques used in
\cite{AAK,BdPS} which crucially exploited structural results
describing type $I$ von Neumann algebras.

\begin{corollary}
\label{c2} Let $\mathcal{M}$ be a $W^*$-algebra and let
$\mathcal{A}$ be a linear subspace in $LS(\mathcal{M})$, such that
$\mathcal{A}^*=\mathcal{A},\ x\in\mathcal{A} \Leftrightarrow
|x|\in\mathcal{A},\ 0\leq x\leq y\in\mathcal{A}\Rightarrow
x\in\mathcal{A}$. Fix $a\in LS(\mathcal{M})$ and consider inner
derivation $\delta=\delta_a$ on the algebra $LS(\mathcal{M})$ given
by $\delta(x)=[a,x]$, $x\in LS(\mathcal{M})$. If
$\delta(\mathcal{M})\subseteq\mathcal{A}$, then there exists an element
$d\in\mathcal{A}$ such that $\delta(x)=[d,x]$.

\end{corollary}
\begin{proof}
Let $a=a_1+i a_2$, where $a_1=\mathfrak{Re}(a)$ and $a_2=\mathfrak{Im}(a)$. We have
$2[a_1,x]=[a+a^*,x]=[a,x]-[a,x^*]^*=\mathcal{A}-\mathcal{A}^*\subseteq\mathcal{A}$
for any $x\in\mathcal{M}$. Analogously, $[a_2,x]\in\mathcal{A}$
for any $x\in\mathcal{M}$. By Theorem \ref{main}, there is a
element $c_k\in Z_h(LS(\mathcal{M}))$ and a unitary element
$u_k\in U(\mathcal{M})$, such that $|[a_k,u_k]|\geq
1/2|a_k-c_k|$ for $k=1,2$. The assumption on $\mathcal{A}$ guarantees that
$a_k-c_k\in\mathcal{A}$, for $k=1,2$. Setting
$d=(a_1-c_1)+i(a_2-c_2)$, we
deduce that $d\in\mathcal{A}$ and $\delta=[d,\cdot]$.
\end{proof}

Consider the following classical example
(see e.g. \cite[Lemmas 3.1 and 3.2]{JP}): $\mathcal{M}=B(H)$ and
$\mathcal{A}$ is the algebra of all compact operators on $H$.
Suppose that an element $a=a^*\in \mathcal{M}$ is such that
$\delta_a:\mathcal{M}\to \mathcal{A}$. Then the result of
Corollary \ref{c2} asserts that there exists $\lambda \in
\mathbb{R}$ such that $a-\lambda \mathbf{1}$ is a compact
operator. An important example extending this
classical result can be obtained as follows. Let a semi-finite von
Neumann  algebra $\mathcal M$ be equipped with a faithful normal
semi-finite trace $\tau$. Let $x\in S(\mathcal{M},\tau)$. The set
$S_0(\mathcal{M},\tau)$ of all $\tau$-compact operators in
$LS(\mathcal{M})$ is defined as the subset of all $x\in
S(\mathcal{M},\tau)$ such that $\lim_{t\to \infty}\mu _{t}(x)=0$
(see the definition the generalized singular value function $\mu$
below). The result of Corollary \ref{c2} asserts, in particular,
that for any $a\in LS(\mathcal{M})$ such that $\delta_a:
LS(\mathcal{M})\to S_0(\mathcal{M},\tau)$ there exists an element
$c\in LS(Z(\mathcal{M}))$ such that $a-c\in
S_0(\mathcal{M},\tau)$.


Numerous examples of absolutely solid
subspaces $\mathcal{A}$ in $LS(\mathcal{M})$ satisfying the
assumptions of the preceding corollary are given by $\mathcal{M}$-bimodules of $LS(\mathcal{M})$.

\begin{definition}
A linear subspace $E$ of $LS(\mathcal{M})$, is
called an $\mathcal{M}$-bimodule of local measurable operators
if $uxv\in E$ whenever $x\in E$ and $u,v\in \mathcal{M}$. If an
$\mathcal{M}$-bimodule $E$ is equipped with a (semi-) norm
$\left\Vert \cdot \right\Vert _{E}$, satisfying
\begin{equation}
\left\Vert uxv\right\Vert _{E}\leq \left\Vert u\right\Vert
_{\mathcal{M}}\left\Vert v\right\Vert _{\mathcal{M}}\left\Vert
x\right\Vert _{E},\ \ \ x\in E,\ u,v\in \mathcal{M}\text{,}
\label{ChIVeq21}
\end{equation}%
then $E$ is called a (semi-) normed $\mathcal{M}$-bimodule of local measurable operators.
\end{definition}

We omit a straightforward verification of the fact that every $\mathcal{M}$-bimodule
of locally measurable operators satisfies the assumption of Corollary \ref{c2}.

The best known examples of normed $\mathcal{M}$-bimodules of $LS\left(\mathcal{M}\right)$ are given by
the so-called symmetric operator spaces (see e.g. \cite{DDP0, SC, KS}).
We briefly recall relevant definitions (for more detailed information we refer to \cite{KS} and references therein).

Let $L_0$ be a space of Lebesgue measurable functions either on
$(0,1)$ or on $(0,\infty)$, or on $\mathbb{N}$ finite almost everywhere
(with identification $m-$a.e.). Here $m$ is Lebesgue measure or else
counting measure on $\mathbb{N}$. Define $S$ as the subset of $L_0$
which consists of all functions $x$ such that $m(\{|x|>s\})$ is finite for some $s.$

Let $E$  be a Banach space of real-valued Lebesgue measurable
functions either on $(0,1)$ or $(0,\infty)$ (with identification
$m-$a.e.) or on $\mathbb{N}$. The space $E$ is said to be {\it
absolutely solid} if $x\in E$ and $|y|\leq |x|$, $y\in L_0$
implies that $y\in E$ and $||y||_E\leq||x||_E.$

The absolutely solid space $E\subseteq S$ is said to be
{\it symmetric} if for every $x\in E$ and every $y$ the assumption
$y^*=x^*$ implies that $y\in E$ and $||y||_E=||x||_E$.

Here, $x^*$ denotes the non-increasing right-continuous rearrangement of $x$ given by
$$x^*(t)=\inf\{s\geq0:\ m(\{|x|\geq s\})\leq t\}.$$
In the case when $x$ is a sequence we denote by $x^*$ the usual
decreasing rearrangement of the sequence $|x|$.

If $E=E(0,1)$ is a symmetric space on $(0,1),$ then
$$L_{\infty}\subseteq E\subseteq L_1.$$

If $E=E(0,\infty)$ is a symmetric space on $(0,\infty),$ then
$$L_1\cap L_{\infty}\subseteq E\subseteq L_1+L_{\infty}.$$

If $E=E(\mathbb{N})$ is a symmetric space on $\mathbb{N},$ then
$$\ell_1\subseteq E\subseteq \ell_{\infty},$$
where $\ell_1$ and $\ell_\infty$ are classical spaces of all
absolutely summable and bounded sequences respectively.

Let a semi-finite von Neumann  algebra $\mathcal M$ be equipped with a
faithful normal semi-finite trace $\tau$. Let $x\in S(\mathcal{M},\tau)$.
The generalized singular value function of $x$ is $\mu (x):t\rightarrow \mu _{t}(x),$ where, for
$0\leq t<\tau (\bf {1})$

\begin{equation*}
\mu _{t}(x)=\inf \{s\geq 0\mid \tau (e^{\left\vert x\right\vert
}(s,\infty )\leq t\}.
\end{equation*}

Consider $\mathcal{M} =L^\infty([0,\infty))$ as an abelian von
Neumann algebra acting via multiplication on the Hilbert space
${\mathcal {H}}=L^2(0,\infty)$, with the trace given by
integration with respect to Lebesgue measure $m$. It is easy to
see that the set of all $\tau$-measurable operators affiliated
with $\mathcal{M}$ consists of all measurable functions on
$[0,\infty)$ which are bounded except on a set of finite measure,
that is $S(\mathcal{M}, \tau)=S$ and that the generalized singular
value function $\mu(f)$ is precisely the decreasing rearrangement
$f^*$.

If $\mathcal{M} ={B}({H})$ (respectively, $\ell_\infty(\mathbb N)$)
and $\tau $ is the standard trace ${\rm Tr}$ (respectively, the counting
measure  on $\mathbb{N}$), then it is not difficult to see
that $S(\mathcal{M}, \tau) =\mathcal{M}$. In this case, for
$x\in S(\mathcal{M}, \tau)$ we have
$$
 \mu _n(x)=\mu _t(x), \quad t\in (n-1,n],\quad  n=0,1,2,\dots .
$$
For $\mathcal M={B}({H})$ the sequence $\{\mu
_n(T)\}_{_{n=1}}^{\infty }$ is just the sequence of singular
values $(s_n(T))_{n=1}^{\infty}.$

\begin{definition}\label{opspace}
Let $\mathcal{E}$ be a linear subset in $S({\mathcal{M}, \tau})$
equipped with a norm $\|\cdot\|_{\mathcal{E}}$. We say that
$\mathcal{E}$ is a \textit{symmetric operator space} (on
$\mathcal{M}$, or in $S({\mathcal{M}, \tau})$) if $x\in \mathcal{E}$ and every $y\in
S({\mathcal{M}, \tau})$ the assumption $\mu(y)\leq \mu(x)$ implies
that $y\in \mathcal{E}$ and $\|y\|_\mathcal{E}\leq
\|x\|_\mathcal{E}$.
\end{definition}
} The fact that every symmetric operator space $\mathcal{E}$ is
(an absolutely solid) $\mathcal{M}$-bimodule of
$S\left(\mathcal{M}, \tau \right)$ is well known (see e.g.
\cite{SC, KS} and references therein). In the special case, when
$\mathcal{M}=B(H)$ and $\tau$ is a standard trace ${\rm Tr}$, the
notion of symmetric operator space introduced in Definition
\ref{opspace} coincides with the notion of symmetric operator
ideal given in  Definition \ref{opideal}.

There exists a strong connection between symmetric function and
operator spaces recently exposed in \cite{KS} (see earlier results
in \cite{Schatten, GK1, GK2, Simon}).

Let $E$ be a symmetric function space on the interval $(0,1)$ (respectively, on the semi-axis or on $\mathbb{N}$)
and let $\mathcal{M}$ be a type $II_1$ (respectively, $II_{\infty}$ or type $I$) von Neumann algebra. Define
$$E(\mathcal{M},\tau):=\{S\in S(\mathcal{M},\tau):\ \mu_t(S)\in E\},\ \|S\|_{E(\mathcal{M},\tau)}:=\|\mu_t(S)\|_E.$$
Main results of \cite{KS} assert that $(E(\mathcal{M},\tau),
\|\cdot\|_{E(\mathcal{M},\tau)})$ is a symmetric operator space.
If $E=L_p$, $1\leq p<\infty$, then $(E(\mathcal{M},\tau),
\|\cdot\|_{E(\mathcal{M},\tau)})$ coincides with the classical
noncommutative $L_p$-space associated with the algebra
$(\mathcal{M},\tau)$. If $\mathcal{M}$ is a semifinite  atomless
von Neumann algebra, then the converse result also holds
\cite{SC}. That is, if $\mathcal{E}$ is a symmetric operator space
on $\mathcal{M}$, then
$$
E(0,\infty):=\{f\in S_0((0,\infty)):\ f^*=\mu(x)\ \text{for some}\ x\in \mathcal{E}\},\ \|f\|_E:=\|x\|_{\mathcal{E}}
$$
is a symmetric function space on $(0,\tau(\textbf{1})$. It is
obvious that $\mathcal{E}=E(\mathcal{M},\tau)$. Similarly, if $E=E(\mathbb{N})$ is a symmetric sequence space on $\mathbb{N},$ and the
algebra $\mathcal {M}$ is a type $I$ factor with standard trace, then (see \cite{KS}) setting
$$\mathcal{E}:=\{S\in \mathcal{M}:\ (s_n(S))_{n=1}^{\infty}\in E\},\ \|S\|_{\mathcal{E}}:=\|(s_n(S))_{n=1}^{\infty}\|_E$$
yields a symmetric operator ideal. Conversely, every symmetric operator ideal $\mathcal{E}$ in $\mathcal {M}$
defines a unique symmetric sequence space $E=E(\mathbb{N})$ by setting
$$
E:=\{a=(a_n)_{n=1}^{\infty}\in \ell_\infty:\ a^*=(s_n(S))_{n=1}^{\infty}\ \text{for some}\ S\in \mathcal{E}\},\ \|a\|_E:=\|S\|_{\mathcal{E}}
$$
We are now fully equipped to provide a full analogue of
Corollaries \ref{c1} and \ref{c11}.

\begin{corollary}
\label{c3} Let $\mathcal{M}$ be a semi-finite $W^*$-algebra and
let $\mathcal{E}$ be a symmetric operator space. Fix $a=a^*\in
S(\mathcal{M})$ and consider inner derivation $\delta=\delta_a$ on
the algebra $LS(\mathcal{M})$ given by $\delta(x)=[a,x]$, $x\in
LS(\mathcal{M})$. If $\delta(\mathcal{M})\subseteq \mathcal{E}$,
then there exists $d\in \mathcal{E}$ satisfying the inequality
$\|d\|_{\mathcal{E}}\leq\|\delta\|_{\mathcal{M}\rightarrow
\mathcal{E}}$ and such that $\delta(x)=[d,x]$.
\end{corollary}

\begin{proof}
The existence of $d\in \mathcal{E}$ such that $\delta(x)=[d,x]$
follows from Corollary \ref{c2}. Now, if $u\in U(\mathcal{M})$,
then $\|\delta(u)\|_{\mathcal{E}}=\|du-ud\|_{\mathcal{E}}\leq
\|du\|_{\mathcal{E}}+\|ud\|_{\mathcal{E}}=2\|d\|_{\mathcal{E}}$.
Hence, if $x\in \mathcal{M}_1=\{x\in\mathcal{M}:\ \|x\|\leq 1\}$,
then $x=\sum_{i=1}^4 \alpha_i u_i$, where $u_i\in U(\mathcal{M})$
and $|\alpha_i|\leq 1$ for $i=1,2,3,4$, and so
$\|\delta(x)\|_{\mathcal{E}}\leq \sum_{i=1}^4\|\delta(\alpha_i
u_i)\|_{\mathcal{E}}\leq 8\|d\|_{\mathcal{E}}$, that is
$\|\delta\|_{\mathcal{M}\rightarrow \mathcal{E}}\leq
8\|d\|_{\mathcal{E}}<\infty$.

The final assertion is established exactly as in the proof of
Corollary \ref{c11}.
\end{proof}

An illustration of the result above complementing the example
given after Corollary \ref{c2},  can be obtained when the space
$E$ is given by the norm closure of the subspace $L_1\cap
L_{\infty}$ in the space $L_1+L_{\infty}$. In this case, the space
$\mathcal{E}=E(\mathcal{M},\tau)$ can be equivalently described as
the set of all $x\in L_1+L_{\infty}(\mathcal{M},\tau)$ such that
$\lim _{t\to \infty} \mu_t(x)=0$. This space is a normed
counterpart of the space $S_0(\mathcal{M},\tau)$ of all
$\tau$-compact operators in $LS(\mathcal {M})$.

\section{The outline of the proof of Theorem
\ref{main}}\label{the outline}

The next two examples demonstrate respectively two main ideas behind the proof of Theorem \ref{main}.  The first example yields the proof for the case when $\mathcal{M}$ coincides with algebra $M_n(\mathbb{C})$ of all $n\times n$ complex matrices.

\begin{example}\label{e1} The assertion of the Theorem \ref{main} holds for the case $\mathcal{M}=M_n(\mathbb{C})$, $n\in \mathbb{N}$
\end{example}
\begin{proof}
Fix
$a=a^*\in M_n(\mathbb{C})$ and select a unitary matrix $v\in M_n(\mathbb{C})$ such that
$$v^*av=\begin{vmatrix}
\lambda_1 & 0&\hdotsfor{2} & 0\\
0 &\lambda_2&\hdotsfor{2} & 0\\
\hdotsfor{5}\\
\hdotsfor{5}\\
0 & 0&\hdotsfor{2} & \lambda_n
\end{vmatrix}\in M_n(\mathbb{C}),
$$
where $\lambda_1\leq\lambda_2\leq...\leq\lambda_n$.

Let the unitary matrix $u\in M_n(\mathbb{C})$ be counter-diagonal,
that is
$$u=\begin{vmatrix}
0 & \hdotsfor{2} &0 & 1\\
0 &\hdotsfor{2} & 1& 0\\
\hdotsfor{5}\\
\hdotsfor{5}\\
1 & 0&\hdotsfor{2} & 0
\end{vmatrix},
$$
and observe that
$$u^*v^*avu=\begin{vmatrix}
\lambda_n & 0&\hdotsfor{2} & 0\\
0 &\lambda_{n-1}&\hdotsfor{2} & 0\\
\hdotsfor{5}\\
\hdotsfor{5}\\
0 & 0&\hdotsfor{2} & \lambda_1
\end{vmatrix}.
$$
Therefore,
$$|[a,vu]|=|u^*v^*avu-a|=\begin{vmatrix}
|\lambda_n-\lambda_1| & 0&\hdotsfor{2} & 0\\
0 &|\lambda_{n-1}-\lambda_2|&\hdotsfor{2} & 0\\
\hdotsfor{5}\\
\hdotsfor{5}\\
0 & 0&\hdotsfor{2} & |\lambda_1-\lambda_n|
\end{vmatrix}.
$$
If $n$ is odd, then for all $1\leq k\leq n$ we have
$$|\lambda_k-\lambda_{n+1-k}|=|\lambda_k-\lambda_0|+|\lambda_{n+1-k}-\lambda_0|$$ for $\lambda_0=\lambda_{(n+1)/2}$.

If $n$ is even, then then for all $1\leq k\leq n$ we have
$$|\lambda_k-\lambda_{n+1-k}|=|\lambda_k-\lambda_0|+|\lambda_{n+1-k}-\lambda_0|$$
for every $\lambda_0\in[\lambda_{n/2},\lambda_{n/2+1}].$

Therefore, for every $n\in \mathbb{N}$, we have
$$|[v^*av,u]|=u^*v^*|a-\lambda_0 \mathbf{1}|vu+v^*|a-\lambda_0 \mathbf{1}|v.$$
Then
$|[a,vuv^*]|=|(vuv^*)^*a(vuv^*)-a|=v|u^*v^*avu-v^*av|v^*=v|[v^*av,u]|v^*=(vuv^*)^*|a-\lambda_0
\mathbf{1}|(vuv^*)+|a-\lambda_0 \mathbf{1}|,\
(vuv^*)^2=\mathbf{1}$.

This completes the proof of the equality \eqref{first main}. The second assertion of Theorem \ref{main} trivially follows from the first (see Remark \ref{trivial}).
\end{proof}

The idea of the proof of the second part of Theorem \ref{main} is demonstrated in Example \ref{e1}.

\begin{example}\label{e2} Let $\mathcal{M}$ be a $\sigma$-finite
factor of type $I_\infty$ or $II_\infty$. Then $\mathcal{M}$
contain a family of pairwise orthogonal and pairwise equivalent
projections $\{p_n\}_{n=1}^\infty$, such that $\sum_{n=1}^\infty
p_n=\mathbf{1}$. Let $\mathbb{R}_+\ni\lambda_n\downarrow 0$. Set
$a=\sum_{n=1}^\infty \lambda_n p_n$ (this series converges in the
strong operator topology). Using the same arguments as in
\cite{BS1}, it is easy to show that there is no
$\lambda_0\in\mathbb{C}$ and $u\in U(\mathcal{M})$ such that
$|[a,u]|\geq |a-\lambda_0\mathbf{1}|$. Nevertheless, the part
(ii) of Theorem \ref{main} holds.
\end{example}

\begin{proof} We refer the reader to \cite{BS1} for the proof of the first assertion and pass to the proof of the second one.

Fix $\varepsilon>0$. For every $n\in\mathbb{N}$ there exist
infinitely many $m\in\mathbb{N}$ such that
$\lambda_m<\varepsilon\lambda_n$. Hence, the set $\mathbb{N}$ can
be split up into the set of all pairs $\{n_k,m_k\}_{k=1}^\infty$
such that $\lambda_{m_k}<\varepsilon\lambda_{n_k}$.
For every $k\ge 1$, select a partial isometry $v_{n_k m_k}$ such
that $v_{n_k m_k}^*v_{n_k m_k}=p_{m_k},\ v_{n_k m_k}v_{n_k
m_k}^*=p_{n_k}$. Clearly, the projections $p_{n_k}, p_{m_k}$ and
the partial isometry $v_{n_k m_k}$ generate a $*$-subalgebra in
$\mathcal{M}$ which is $*$-isomorphic to $M_2(\mathbb{C})$.
Without loss of generality, under this $*$-isomorphism, we equate
$$
p_{n_k}=\begin{vmatrix}1&0\\0&0\end{vmatrix},
p_{m_k}=\begin{vmatrix}0&0\\0&1\end{vmatrix}, v_{n_k
m_k}=\begin{vmatrix}0&1\\0&0\end{vmatrix}.
$$
Next, consider a $W^*$-subalgebra of $\mathcal{M}$ generated by the elements $p_{n_k},
p_{m_k},v_{n_k m_k}$, $k\ge 1$ which we identify with
$\bigoplus_{k=1}^{\infty}M_2(\mathbb{C})$. We have
$$a=\bigoplus_{k=1}^{\infty}\begin{vmatrix}\lambda_{n_k}&0\\0&\lambda_{m_k}\end{vmatrix}.$$
Set
$$u=\bigoplus_{k=1}^{\infty}\begin{vmatrix}0&1\\1&0\end{vmatrix}.$$
It clearly follows that
$$|[u,a]|=\bigoplus_{k=1}^{\infty}\begin{vmatrix}|\lambda_{n_k}-\lambda_{m_k}|&0\\0&|\lambda_{n_k}-\lambda_{m_k}|\end{vmatrix}\geq(1-\varepsilon)\bigoplus_{k=1}^{\infty}\begin{vmatrix}|\lambda_{n_k}|&0\\0&|\lambda_{m_k}|\end{vmatrix}=(1-\varepsilon)|a|.$$
\end{proof}


The idea of the  proof of the
Theorem \ref{main} consists in the splitting of the identity in
$\mathcal{M}$ into the sum of three pairwise orthogonal central
projections $p_0,p_-,p_+$ with certain properties. For the reduced
algebra $\mathcal{M}p_0$ the method from example \ref{e1} will be
applied,  and for the reduced algebras $\mathcal{M}p_-$  and
$\mathcal{M}p_+$ the method from example \ref{e2} will be
adjusted.

\section{The proof of Theorem \ref{main}}\label{the proof}

For better readability, we break the proof into the following series of
lemmas.

Until the end of  the proof, we fix an arbitrary element $a\in
LS_h(\mathcal{M})$.

For projections $p,q\in P(\mathcal{M})$ we assume $p\prec\prec q$,
if $pz\prec qz$ for every $0<z\leq \mathbf{c}(p)\vee
\mathbf{c}(q),\ z\in P(Z(\mathcal{M}))$. Set $p\succ\succ q$,
if $q\prec\prec p$.

We recall the following comparison result

\begin{theorem}
\label{t11} (\cite[Theorem 6.2.7]{KR2})
Let $e$ and $f$ be projections in a von Neumann algebra $\mathcal{M}$.
There are unique orthogonal central projections $p$ and $q$ maximal with
respect to the properties $qe \sim qf$, and, if $p_0$ is a non-zero
central subprojection of $p$, then $p_0 e \prec p_0 f$. If $r_0$ is a non-zero
central subprojection of $\mathbf{1}-p-q$, then $r_0 f \prec r_0 e$.
\end{theorem}

The following form of the preceding theorem will be more convenient for our purposes.

\begin{theorem}
\label{t2} Let $\mathcal{M}$ be a $W^*$-algebra and $p,q\in
P(\mathcal{M})$. Then there exists a unique triple of pairwise
orthogonal projections  $z_-,z_0,z_+\in P(Z(\mathcal{M}))$,
such that $z_-+z_0+z_+=\mathbf{1},\ z_-p\prec\prec z_-q,\
z_0p\sim z_0q,\ z_+p\succ\succ z_+q$.
\end{theorem}

\begin{corollary}
\label{c1_0} Let $\mathcal{M}$ be a  $W^*$-algebra and
$p_1,p_2,...,p_n\in P(\mathcal{M})$. Then there exists a family
$\{z_{\sigma}\}_{\sigma \in S_n}$ (here, $S_n$ is the permutation group of
 $n$ elements) of pairwise orthogonal projections from
$P(Z(\mathcal{M}))$ such that $\sum_{\sigma\in
S_n}z_{\sigma}=\mathbf{1}$ и $z_{\sigma} p_{\sigma(1)}\preceq
z_{\sigma} p_{\sigma(2)}\preceq ... \preceq z_{\sigma}
p_{\sigma(n)}$ for every $\sigma \in S_n$.
\end{corollary}

\begin{proof} For every pair of projections $p_i$,
$p_j$, $i<j$ we denote by $z^1_{ij}$ the largest central
projection such that $z^1_{ij}p_i\preceq z^1_{ij}p_j$. Then for
$z^2_{ij}:=\mathbf{1}-z^1_{ij}$ we have $z^2_{ij}p_i\succ\succ
z^2_{ij}p_j$. Let $\mathfrak{B}$ be a Boolean algebra generated by
all central projections $z^1_{ij}$ and  $z^2_{ij}$, $1\leq i<j\leq
n$. Every atom
$z\in\mathfrak{B}$ may be uniquely written as
\begin{equation}\label{atom}
z=\prod_{i=1}^{n-1}\prod_{j=i+1}^n z^{\varepsilon_{ij}}_{ij},
\end{equation}
where $\varepsilon_{ij}\in\{1,2\}$. Observe that the sum of all
elements as in \eqref{atom} is equal to $\mathbf{1}$.

Fix an atom $z\in\mathfrak{B}$ and define on the set
$\{1,...,n\}$ a (linear) order $\leq_z$. To this end, we set for any $i,j,\
i\neq j$, $$i\leq_z j \Leftrightarrow zp_i\preceq z p_j.$$

This definition is correct, since we have $z\leq
z^{\varepsilon_{ij}}_{ij}$ for $i<j$ or $z\leq
z^{\varepsilon_{ji}}_{ji}$ for $i>j$. So, we have
$zp_i\preceq z p_j$ or $zp_i\succeq z p_j$.

Let $\{i_1,...,i_n\}=\{1,...,n\}$, $i_1\leq_z
i_2\leq_z...\leq_z i_n$ and set $\sigma(k)=i_k$, $\sigma\in
S_n$. Then $z=z_\sigma$.
\end{proof}

Let $\Lambda$ be a lattice. Subset $I\subset\Lambda$
is called an $\vee$-ideal ($\wedge$-ideal) in $\Lambda$ if for
every $s,t\in I$ and $u\in\Lambda$, $u\leq s$ ($u\geq s$)
we have $u\in I$ and $s\vee t \in I$ ($s\wedge t \in I$)
\cite{Bir}.

\begin{lemma}
\label{l16} Let $\nabla$ be a complete boolean algebra, let $I$
be a non-zero $\vee$-ideal in $\nabla$.  Then there exists a family
$\{s_\alpha\}_{\alpha\in \Omega}\subset I$ of pairwise disjoint
non-zero elements such that $\bigvee_{\alpha\in
\Omega}s_\alpha=\bigvee I$.
\end{lemma}
\begin{proof}
It follows from the Zorn's lemma that among all families of pairwise disjoint
non-zero elements from $I$ there exists a maximal one with respect to the inclusion.
Let $\{s_\alpha\}_{\alpha\in \Omega}$ be one of these maximal families. It is clear that
$s:=\bigvee_{\alpha\in \Omega}s_\alpha\leq\bigvee I=t$.
Suppose $s<t$. Then there exists $v\in I$ such that
$v\nleqslant s$. In this case $v\wedge s'\neq 0$
(by $s'$ we denote the complement to $s$). Then
$\{v\wedge s'\}\cup \{s_\alpha\}_{\alpha\in \Omega}$
is a family of pairwise disjoint
non-zero elements from $I$. Thus, we have found a contradiction to the assumption that $\{s_\alpha\}_{\alpha\in \Omega}$ is maximal.
So, our assumption $s<t$ fails and we have $s=t$.
\end{proof}

The algebra (over reals) of all self-adjoint elements from the
center of the algebra $LS(\mathcal{M})$ will be denoted by
$Z_h(LS(\mathcal{M}))$. The latter algebra is a lattice with
respect to the order induced from  $LS(\mathcal{M})$. As we have
already explained in the Preliminaries the $*$-algebra
$Z(LS(\mathcal{M}))$ is $*$-isomorphic to the $*$-algebra
$L^0(\Omega,\Sigma,\mu)$ of all measurable complex valued
functions on a measurable space with a locally finite measure.
This isomorphism yields the assertions of Lemmas \ref{l17} and
\ref{l14_1} below.

\begin{lemma}
\label{l17} Let $\Omega:=\{c_j\}_{j\in J}$ be a family of pairwise disjoint
elements from $Z_h(LS(\mathcal{M}))$. Then there exists an element $c\in Z_h(LS(\mathcal{M}))$ such that
$c\mathbf{s}(c_j)=c_j$ for every $j\in J$.
\end{lemma}

\begin{lemma}
\label{l14_1}
A lattice $Z_h(LS(\mathcal{M}))$ is conditionally complete.
\end{lemma}

\begin{lemma}
\label{l00} Let $p,q,r\in P(\mathcal{M})$, $p<q$, $p\prec
r\prec\prec q$. Then there exists an element $r_1\in
P(\mathcal{M})$ such that $r_1\sim r$ and $p<r_1<q$.
\end{lemma}
\begin{proof}  Due to the assumptions, there exists an element
$p_1\in P(\mathcal{M})$ such that $p\sim p_1<r$. Suppose that
$(r-p_1)z\succeq (q-p)z$ for some $0\neq z\in P(Z(\mathcal{M}))$,
$z\leq \mathbf{c}(q)$. Then $rz=(r-p_1)z+p_1 z\succeq
(q-p)z+pz=qz$. This fact contradicts the assumption $r\prec\prec
q$, which yields the estimate $rz\prec qz$ for every $0<z\leq
\mathbf{c}(q)=\mathbf{c}(r)\vee \mathbf{c}(q)$ (the previous
equality follows from the implication $r\prec q \Rightarrow
\mathbf{c}(r)\leq \mathbf{c}(q)$). Hence, from Theorem \ref{t2} we
have that $r-p_1\prec q-p$. So, there exists an element $0\neq p_2\in
P(\mathcal{M})$, such that
 $r-p_1\sim p_2<q-p$. Then $p<p+p_2<q$ and $p+p_2\sim
p_1+(r-p_1)=r$, where the equivalence $p+p_2\sim p_1+(r-p_1)$ is
guaranteed by the implication
$$p\sim p_1,\ p_2\sim r-p_1,\ pp_2=0,\
p_1(r-p_1)=0 \Rightarrow p+p_2\sim p_1+(r-p_1).$$
We are done, by letting $r_1:=p+p_2$.
\end{proof}

\begin{lemma}
\label{l0} Let $p$ be a  properly infinite projection in
$\mathcal{M}$. Then we have

(i). If $P(\mathcal{M}) \ni q_1,...,q_n,... \preceq p$, $q_nq_m=0$
for $n\neq m$, then $\bigvee_{n=1}^\infty q_n \preceq p$.

(ii). If $P(\mathcal{M}) \ni q_1,...,q_n\prec\prec p$, $q_iq_j=0$
for $i\neq j$, then $\bigvee_{i=1}^n q_i \prec\prec p$.

(iii). If $p\succeq \mathbf{1}-p$, then $p\sim \mathbf{1}$.

(iv). If $P(\mathcal{M}) \ni q\prec\prec p,\ qp=pq$, then
$p(\mathbf{1}-q)\sim p$.

(v). If $q\in P(\mathcal{M})$ is a finite projection,
$P(\mathcal{M}) \ni p_n \uparrow p$ and $z_n \in
P(Z(\mathcal{M}))$ are such that $(\mathbf{1}-z_n)p_n\preceq
(\mathbf{1}-z_n)q$ for all $n\in \mathbb{N}$, then
$\bigvee_{n=1}^\infty z_n \geq \mathbf{c}(p)$.

(vi). If $P(\mathcal{M}) \ni q_n \uparrow q\succ\succ p$, $z_n \in
P(Z(\mathcal{M}))$ and $(\mathbf{1}-z_n)p\succeq
(\mathbf{1}-z_n)q_n$ for all $n\in \mathbb{N}$, then
$\bigvee_{n=1}^\infty z_n \geq \mathbf{c}(p)$.
\end{lemma}
\begin{proof}
(i). Since $p$ is a properly infinite projection then there are
pairwise disjoint projections $p_1,...,p_n,...\in P(\mathcal{M})$
such that $p=\bigvee_{n=1}^\infty p_n$, $p_n\sim p$ for all $n\in
\mathbb{N}$. Then $p_n\succeq q_n$ for every $n\in \mathbb{N}$.
Hence, $p=\bigvee_{n=1}^\infty p_n\succeq \bigvee_{n=1}^\infty
q_n$.

(ii). It follows from Corollary \ref{c1_0} that there exists a
family $\{z_{\sigma}\}_{\sigma \in S_n}$ of pairwise orthogonal
projections in $P(Z(\mathcal{M}))$ such that $\sum_{\sigma\in
S_n}z_{\sigma}=\mathbf{c}(p)$, $z_{\sigma} q_{\sigma(1)}\preceq
z_{\sigma} q_{\sigma(2)}\preceq ... \preceq z_{\sigma}
q_{\sigma(n)}$ for all $\sigma \in S_n$.
We shall consider only those elements $\sigma \in S_n$, for which
$z_{\sigma}\neq 0$. If $z_{\sigma} q_{\sigma(n)}$ is a finite
projection, then $\bigvee_{j=1}^n z_{\sigma} q_{\sigma(j)}$ is
also finite and $\bigvee_{j=1}^n z_{\sigma} q_{\sigma(j)}
\prec\prec p z_{\sigma}$. Indeed, it follows from (i) that
$\bigvee_{j=1}^n z_{\sigma} q_{\sigma(j)} \preceq p z_{\sigma}$.
Moreover, we have $p z_{\sigma}$ is a properly infinite
projection. So, there does not exist central projection
 $0<z\leq z_{\sigma},$ such that$\bigvee_{j=1}^n z
q_{\sigma(j)}\neq 0$ that $\bigvee_{j=1}^n z q_{\sigma(j)} \sim p
z$ since, in this case, we would have that $0\neq p z$ is finite. If $z_{\sigma}
q_{\sigma(n)}$ is a properly infinite projection than, according to
(i), $\bigvee_{j=1}^n z_{\sigma} q_{\sigma(j)} \preceq z_{\sigma}
q_{\sigma(n)} \prec\prec p z_{\sigma}$. Hence, $\bigvee_{i=1}^n
q_i=\bigvee_{i=1}^n \sum_{\sigma\in S_n} z_{\sigma} q_i =
\sum_{\sigma\in S_n}\bigvee_{j=1}^n z_{\sigma} q_{\sigma(j)}
\prec\prec \sum_{\sigma\in S_n}p z_{\sigma} =p$.

(iii). Since, $\mathbf{1}=p+(\mathbf{1}-p)$, the estimate
$\mathbf{1}\preceq p$ follows from (i).

(iv). We have $p=qp+(\mathbf{1}-q)p,\ qp\leq q \prec\prec p$.
Suppose $(\mathbf{1}-q)pz\prec\prec pz$ for some $0\neq z\in
P(Z(\mathcal{M}))$. Then, according to (ii), $p\prec\prec p$,
which is false. Hence, $(\mathbf{1}-q)p\succeq p$. So, since
$(\mathbf{1}-q)p\leq p$, we conclude $(\mathbf{1}-q)p\sim p$.

(v). Let $z_0:=\bigwedge_{n=1}^\infty (\mathbf{1}-z_n)$. Then,
it is clear, $z_0\leq \mathbf{1}-z_n$. We shall now use a well-known implication
\begin{equation}
\label{well-known}
e\preceq f,\ z\in
P(Z(\mathcal{M}))\Rightarrow ze\preceq zf.
\end{equation}
By the assumption, we have $(\mathbf{1}-z_n)p_n\preceq
(\mathbf{1}-z_n)q$ and therefore, by \eqref{well-known}
$z_0p_n=z_0(\mathbf{1}-z_n)p_n\preceq z_0(\mathbf{1}-z_n)q=z_0q$
for all $n\in\mathbb{N}$. Then all projections $z_0p_n$ are finite
and it follows from (\cite{Tak}, Chapter V, Lemma 2.2) that
$z_0p\preceq z_0q$. In this case $z_0p$ is a finite projection.
Since $p$ is a properly infinite projection, we conclude $z_0p=0$.
The latter trivially implies that $p(\mathbf{1}-z_0)=p$ and since
$\mathbf{c}(p)\leq \mathbf{1}-z_0$, we also obtain $z_0
\mathbf{c}(p)=0$. So,
$$(\mathbf{1}-\mathbf{c}(p))\vee\bigvee_{n=1}^\infty
z_n=(\mathbf{1}-\mathbf{c}(p))\vee(\mathbf{1}-z_0)=\mathbf{1},$$
in particular, $\bigvee_{n=1}^\infty z_n\geq \mathbf{c}(p)$.

(vi). Let $z_0:=\bigwedge_{n=1}^\infty (\mathbf{1}-z_n)$. Then,
from the conditions $\mathbf{1}-z_n\geq z_0,\
(\mathbf{1}-z_n)p\succeq (\mathbf{1}-z_n)q_n$ we have
 $z_0p\succeq z_0q_n$ for all
$n\in \mathbb{N}$. So, it follows from (i) that  $z_0p=0$ or
$z_0p\succeq z_0q_1+\bigvee_{n=1}^\infty
(z_0q_{n+1}-z_0q_n)=z_0q$. However, due to the assumption $q\succ\succ p$,
the estimate $z_0p\succeq z_0q$ can hold only when
$z_0\mathbf{c}(q)=0$. Since
$\mathbf{c}(p)\leq \mathbf{c}(q)$, we have
$z_0\mathbf{c}(p)=0$ in any case. Hence, $\bigvee_{n=1}^\infty z_n =
\mathbf{1}-z_0\geq \mathbf{c}(p)$.
\end{proof}

Let $c\in Z_h(LS(\mathcal{M}))$. We set

$$e^a_z(-\infty,c):=\mathbf{s}((c-a)_+),\ e^a_z(c,+\infty):=\mathbf{s}((a-c)_+),$$
$$e^a_z(-\infty,c]:=\mathbf{1}-e^a_z(c,+\infty),\ e^a_z[c,+\infty):=\mathbf{1}-e^a_z(-\infty,c).
$$

Observe that all the projections defined above belong to a commutative $W^*$-subalgebra of $\mathcal{M}$, generated by $Z(\mathcal{M})$ and
spectral projections of the element $a$. Having this observation in mind, for any
$c_1,c_2\in Z_h(LS(\mathcal{M}))$, such that $c_1\leq c_2$, we set

$$e^a_z[c_1,c_2):=e^a_z(-\infty,c_2)e^a_z[c_1,+\infty),\ e^a_z(c_1,c_2]:=e^a_z(-\infty,c_2]e^a_z(c_1,+\infty),$$
$$e^a_z[c_1,c_2]:=e^a_z(-\infty,c_2]e^a_z[c_1,+\infty),\ e^a_z(c_1,c_2):=e^a_z(-\infty,c_2)e^a_z(c_1,+\infty).$$

Finally, we set $$e^a_z\{c\}=e^a_z[c,c].$$

Observe that our spectral measure is analogous to the construction
given in \cite[Definition 2.4]{Hl} (for the case of
unbounded locally measurable self-adjoint operator).


\begin{lemma}
\label{l11}
Fix $c\in Z_h(LS(\mathcal{M}))$.

(i). If $c_1,c_2\in Z_h(LS(\mathcal{M}))$, $c_1\leq c_2$ then $e^a_z(-\infty,c_1)\leq e^a_z(-\infty,c_2)$ and
$e^a_z(c_1,+\infty)\geq e^a_z(c_2,+\infty)$;

(ii). If $z\in P(Z(\mathcal{M}))$ then
$e^a_z(-\infty,c)z=e^a_z(-\infty,cz)z=e^{az}_z(-\infty,cz)z$ and
$e^a_z(c,+\infty)z=e^a_z(cz,+\infty)z=e^{az}_z(cz,+\infty)z$;

(iii). $a e^a_z(-\infty,c]\leq ce^a_z(-\infty,c]$;

(iv). $ae^a_z[c,+\infty)\geq ce^a_z[c,+\infty)$;

(v). $ae^a_z\{c\}=ce^a_z\{c\}$;

(vi). If $\{c_\alpha\}_{\alpha\in I}\subset Z_h(LS(\mathcal{M}))$ and $c=\bigvee_{\alpha\in I} c_\alpha$ then $\bigvee_{\alpha\in I} e^a_z(-\infty,c_\alpha) = e^a_z(-\infty,c)$;

(vii). If $\{c_\alpha\}_{\alpha\in I}\subset Z_h(LS(\mathcal{M}))$ and $c=\bigwedge_{\alpha\in I} c_\alpha$ then $\bigvee_{\alpha\in I} e^a_z(c_\alpha,+\infty) = e^a_z(c,+\infty)$.
\end{lemma}

\begin{proof}
(i). Since $c_1-a \leq c_2-a$, then $\mathbf{s}((c_1-a)_+)\leq
\mathbf{s}((c_2-a)_+)$ and $\mathbf{s}((a-c_1)_+)\geq
\mathbf{s}((a-c_2)_+)$.

(ii).
$e^a_z(-\infty,c)z=\mathbf{s}((c-a)_+)z=\mathbf{s}((cz-a)_+)z=\mathbf{s}((cz-az)_+)z$.
The second set of equalities is proven in the same way.

(iii).
$(c-a)e^a_z(-\infty,c]=(c-a)(\mathbf{1}-e^a_z(c,+\infty))=(c-a)+(a-c)\mathbf{s}((a-c)_+)=(c-a)+(a-c)_+
= (a-c)_-\geq 0$.

The proof of (iv) is the same.

(v). From (iii) and (iv) we have that $ae^a_z\{c\}\leq ce^a_z\{c\}$ и $ae^a_z\{c\}\geq ce^a_z\{c\}$. Hence, $ae^a_z\{c\}=ce^a_z\{c\}$.

(vi). Since, $(\mathbf{1}-e^a_z(-\infty,c))(c_\alpha-a)\leq
(\mathbf{1}-e^a_z(-\infty,c))(c-a)\leq 0$ then
$\mathbf{1}-e^a_z(-\infty,c)\leq
\mathbf{1}-e^a_z(-\infty,c_\alpha)$ or
$e^a_z(-\infty,c_\alpha)\leq e^a_z(-\infty,c)$ for every $\alpha\in
I$. Let $q=e^a_z(-\infty,c)- \bigvee_{\alpha\in I}
e^a_z(-\infty,c_\alpha)$. Then $(c-a)q\geq 0$. On the other hand,
$(c_\alpha-a)q\leq 0$ for every $\alpha\in I$. Hence,
$(c-a)q=\bigvee_{\alpha\in I} (c_\alpha-a)q\leq 0$. Then
$(c-a)q=0$. So, $q=0$ since $q\leq \mathbf{s}((c-a)_+)$ and $q$
commute with $c-a$.

(vii). Since $c\leq c_\alpha$, then $e^a_z(c_\alpha,+\infty)\leq e^a_z(c,+\infty)$ for every $\alpha\in I$ (see the beginning of the proof of (vi)).
Let $q=e^a_z(c,+\infty)- \bigvee_{\alpha\in I} e^a_z(c_\alpha,+\infty)$.
Then $q(a-c)\geq 0$. On the other hand, $q(a-c_\alpha)\leq 0$ for every $\alpha\in I$.
So, $q(a-c)=\bigvee_{\alpha\in I}q(a-c_\alpha)\leq 0$. Hence, $q(a-c)=0$ and $q=0$.
\end{proof}

\begin{lemma}
\label{l12}
There exists an element $c\in Z_h(LS(\mathcal{M}))$
such that
\begin{equation}
\label{l100}
pe^a_z(-\infty,c)\prec pe^a_z(c,+\infty), \ \forall p\in P(Z(\mathcal{M})),\ p>0.
\end{equation}
\end{lemma}
\begin{proof}
Since $a\in LS_h(\mathcal{M})$, there exists a set of
pairwise disjoint projections $\{p_n\}_{n=1}^\infty$ from
$P(Z(\mathcal{M}))$ such that $\bigvee_{n=1}^\infty p_n =\mathbf{1}$ and
$ap_n\in S_h(\mathcal{M})$ for every $n\in\mathbb{N}$. Without loss of generality, we mayassume $a\in S_h(\mathcal{M})$.
Indeed, if for every $ap_n$, $n\ge 1$ there exists an element $c^{(n)}\in Z_h(LS(\mathcal{M}p_n))$, satisfying (\ref{l100}),
then, by Lemma \ref{l17}, there exists an element $c\in Z_h(LS(\mathcal{M}))$ satisfying $cp_n=c^{(n)}p_n$, $n\ge 1$ and (\ref{l100}).
So, we may (and shall) consider only the case $a\in S_h(\mathcal{M})$.
The latter assumption guarantees that there exists a scalar $\lambda_1\in\mathbb{R}$ such that
$e^a(-\infty,\lambda_1)$ is a finite projection. Let
$\{\lambda_n\}_{n=1}^\infty \subset \mathbb{R},\
\lambda_n\downarrow -\infty$. By the spectral theorem, we have
$e^a(-\infty,\lambda_n)\downarrow 0$. For every $n\in
\mathbb{N}$ set
$$q_n:=\bigvee\{r\in P(Z(\mathcal{M})):\
re^a(-\infty,\lambda_n)\succeq re^a(\lambda_n,+\infty)\},\
q:=\bigwedge_{n=1}^\infty q_n.$$ Fix $n$ and let $k$
tend to infinity. Then $qe^a(-\infty,\lambda_{n+k})\succeq
qe^a(\lambda_{n+k},+\infty)\geq qe^a(\lambda_n,+\infty)$. Since
all projections $qe^a(-\infty,\lambda_{n+k})$ are finite and
$qe^a(-\infty,\lambda_{n+k})\downarrow 0$, we have
$qe^a(\lambda_n,+\infty)=0$ for every $n\in\mathbb{N}$
(\cite{BdPS}, Lemma 6.11). On the other hand,
$e^a(\lambda_n,+\infty)\uparrow\mathbf{1}$. So,
$qe^a(\lambda_n,+\infty)\uparrow q$. Hence, $q=0$. Then
$\bigvee_{n=1}^\infty (\mathbf{1}-q_n) = \mathbf{1}$.

Let
$$r_1:=\mathbf{1}-q_1,\
r_{n+1}:=\bigvee_{k=1}^{n+1}(\mathbf{1}-q_k)-\bigvee_{k=1}^n(\mathbf{1}-q_k), \ n \ge 1.$$
Then projections $r_n$, $n \geq 1$ are pairwise disjoint, $\bigvee_{n=1}^\infty
r_n=\mathbf{1}$, $r_n\leq \mathbf{1}-q_n$ for all
$n\in\mathbb{N}$. It follows from Lemma \ref{l17} that there exists an element
$c\in Z_h(LS(\mathcal{M}))$ such that $cr_n=\lambda_n r_n$.

Let $0\neq p\in P(Z(\mathcal{M}))$. It follows from Theorem \ref{t2}
that $p=e_1+e_2$, where $e_1,e_2\in P(Z(\mathcal{M}))$, $e_1
e^a_z(-\infty,c)\prec e_1 e^a_z(c,+\infty)$, $e_2
e^a_z(-\infty,c)\succeq e_2 e^a_z(c,+\infty)$. In this case, using Lemma \ref{l11}(ii),
we obtain
\begin{align*}
e_2 r_n e^a(-\infty,\lambda_n)&= e_2 r_n e^a(-\infty,\lambda_n r_n) = e_2 r_n e^a_z(-\infty,cr_n)\\
&=e_2 r_n e^a_z(-\infty,c)\succeq e_2 r_n e^a_z(c,+\infty) \\
&=e_2 r_n e^a_z(cr_n,+\infty) = e_2 r_n e^a(\lambda_nr_n,+\infty) \\
&= e_2 r_n e^a(\lambda_n,+\infty)
\end{align*}
for any $n\in\mathbb{N}$. Due to the definition of the projection $q_n$, we obtain $e_2 r_n \leq
q_n$.  Hence, $e_2 r_n = 0$
(since $e_2 r_n\leq r_n\leq \mathbf{1}-q_n$) for every
$n\in\mathbb{N}$. So, $e_2=\bigvee_{n=1}^\infty e_2 r_n =0$.
Hence, $p=e_1$ and $p e^a_z(-\infty,c)\prec p e^a_z(c,+\infty)$.
\end{proof}

Set
$$\Lambda_-:=\{c\in Z_h(LS(\mathcal{M})):\
pe^a_z(-\infty,c)\prec pe^a_z(c,+\infty),\ \forall p\in P(Z(\mathcal{M})),\ p>0\},
$$
$$\Lambda_+:=\{c\in
Z_h(LS(\mathcal{M})):\ pe^a_z(-\infty,c)\succ
pe^a_z(c,+\infty),\ \forall p\in P(Z(\mathcal{M})),\ p>0\}.$$
It follows from Lemma \ref{l12} that
$\Lambda_-\neq\emptyset$. The proof that $\Lambda_+\neq\emptyset$ is similar.

\begin{lemma}
\label{l111} $\mathbf{c}(e^a_z(c,+\infty))=\mathbf{1}$ for every
$c\in\Lambda_-$ and $\mathbf{c}(e^a_z(-\infty,c))=\mathbf{1}$ for
every $c\in\Lambda_+$.
\end{lemma}
\begin{proof}
Let $c\in\Lambda_-$. Suppose
$p=\mathbf{1}-\mathbf{c}(e^a_z(c,+\infty))>0$. Then
$0=pe^a_z(c,+\infty)\succ pe^a_z(-\infty,c)\geq 0$ which means $0>0$.
This contradiction shows $p=0$.

The proof of the second part is similar.
\end{proof}

\begin{lemma}
\label{l13}
$\Lambda_-$ is a $\vee$-ideal and  $\Lambda_+$ is a $\wedge$-ideal in the lattice $Z_h(LS(\mathcal{M}))$.
\end{lemma}
\begin{proof}
We will show only that $\Lambda_-$ is a $\vee$-ideal in  $Z_h(LS(\mathcal{M}))$. The proof of the second part is analogous.

Let $c_1,c_2\in\Lambda_-$.
If $c_1 \geq c\in Z_h(LS(\mathcal{M}))$
then from Lemma \ref{l11} (i) and the definition of $\Lambda_-$ we have $pe^a_z(-\infty,c)\leq pe^a_z(-\infty,c_1)\prec
pe^a_z(c_1,+\infty)\leq p e^a_z(c,+\infty)$
for every $0\neq p \in P(Z(\mathcal{M}))$. So, $c\in\Lambda_-$.

We set $q:=\mathbf{s}((c_2-c_1)_+)$ and observe that $q\in
P(Z(\mathcal{M}))$, and that $qc_1\leq qc_2,\
(\mathbf{1}-q)c_1\geq (\mathbf{1}-q)c_2$. Then $c_1\vee
c_2=(\mathbf{1}-q)c_1+qc_2$. Let $0\neq p\in P(Z(\mathcal{M}))$.
Then, from Lemma \ref{l11} (ii) we have $pe^a_z(-\infty,c_1\vee
c_2)=p(\mathbf{1}-q)e^a_z(-\infty,c_1)+pqe^a_z(-\infty,c_2)$. In
this case $p(\mathbf{1}-q)=0$ or
$p(\mathbf{1}-q)e^a_z(-\infty,c_1)\prec
p(\mathbf{1}-q)e^a_z(c_1,+\infty)$. So, $pq=0$ or
$pqe^a_z(-\infty,c_2)\prec pqe^a_z(c_2,+\infty)$. Hence, in both cases we have
$pe^a_z(-\infty,c_1\vee c_2)=
p(\mathbf{1}-q)e^a_z(-\infty,c_1)+pqe^a_z(-\infty,c_2) \prec
p(\mathbf{1}-q)e^a_z(c_1,+\infty)+pqe^a_z(c_2,+\infty) =
pe^a_z(c_1\vee c_2,+\infty)$. Hence, $c_1\vee c_2 \in \Lambda_-$.
\end{proof}

\begin{lemma}
\label{l14}
Let $c_1\in\Lambda_-,\ c_2\in\Lambda_+$. Then $pc_1<pc_2$ for every $0\neq p\in P(Z(\mathcal{M}))$.
\end{lemma}
\begin{proof}
Suppose a contrary. Then $pc_1\geq pc_2$
for some $0\neq p\in P(Z(\mathcal{M}))$. In this case,
\begin{align*}
pe^a_z(-\infty,c_1)&= p e^a_z(-\infty,pc_1) \quad\text{(by Lemma \ref{l11} (ii) )}\\
&\geq pe^a_z(-\infty,pc_2) \quad\text{(by Lemma \ref{l11} (i) )}\\
&  =p e^a_z(-\infty,c_2)\quad\text{(by Lemma \ref{l11} (i) )}\\
& \succ pe^a_z(c_2,+\infty)\quad\text{(by the definition)}\\
& =pe^a_z(pc_2,+\infty)\quad\text{(by Lemma \ref{l11} (ii) )}\\
& \geq pe^a_z(pc_1,+\infty)\quad\text{(by Lemma \ref{l11} (i) )}\\
&= pe^a_z(c_1,+\infty), \quad\text{(by Lemma \ref{l11} (ii) )}
\end{align*}
which is not true since $c_1\in\Lambda_-$.
\end{proof}

Due to Lemmas \ref{l14_1} and \ref{l14}, we can now set
$$c_0:=\bigvee\Lambda_-\in Z_h(LS(\mathcal{M}))$$
and apply again Lemma \ref{l14} to infer that $c_0\leq c$ for any
$c\in\Lambda_+$.

\begin{lemma}
\label{l15_0}
$c_0-\varepsilon \mathbf{1}\in\Lambda_-$ for every $\varepsilon>0$.
\end{lemma}
\begin{proof}
Suppose a contrary. Then there exists a projection
$0<p\in P(Z(\mathcal{M}))$ such that
$pe^a_z(-\infty,c_0-\varepsilon \mathbf{1})\succeq
pe^a_z(c_0-\varepsilon \mathbf{1},+\infty)$. Choose an arbitrary
element $c\in\Lambda_-$. We shall show that $pc\leq p(c_0-\varepsilon
\mathbf{1})$. If it is not the case, then
$P(Z(\mathcal{M}))\ni q:=\mathbf{s}((p(c_0-\varepsilon
\mathbf{1})-pc)_-)>0,\ q\leq p$ i $qc\geq q(c_0-\varepsilon
\mathbf{1})$. Then
\begin{align*}
 qe^a_z(-\infty,c_0-\varepsilon\mathbf{1})=qe^a_z(-\infty,q(c_0-\varepsilon \mathbf{1})) \quad\text{(by Lemma \ref{l11} (ii) )}\\
 \leq qe^a_z(-\infty,qc) \quad\text{(by Lemma \ref{l11} (i) )}\\
 =qe^a_z(-\infty,c) \quad\text{(by Lemma \ref{l11} (ii) )}\\
 \prec qe^a_z(c,+\infty)\quad\text{(by the definition)}\\
=qe^a_z(qc,+\infty) \quad\text{(by Lemma \ref{l11} (ii) )}\\
\leq qe^a_z(q(c_0-\varepsilon\mathbf{1}),+\infty) \quad\text{(by Lemma \ref{l11} (i) )}\\
=qe^a_z(c_0-\varepsilon \mathbf{1},+\infty) \quad\text{(by Lemma \ref{l11} (ii) )}\\
\end{align*}
On the other hand,
$qe^a_z(-\infty,c_0-\varepsilon \mathbf{1})\succeq
qe^a_z(c_0-\varepsilon \mathbf{1},+\infty)$, since
$pe^a_z(-\infty,c_0-\varepsilon \mathbf{1})\succeq
pe^a_z(c_0-\varepsilon \mathbf{1},+\infty)$ и
$P(Z(\mathcal{M}))\ni q\leq p$. However the inequalities $qe^a_z(-\infty,c_0-\varepsilon
\mathbf{1})\prec qe^a_z(c_0-\varepsilon \mathbf{1},+\infty)$ and
$qe^a_z(-\infty,c_0-\varepsilon \mathbf{1})\succeq
qe^a_z(c_0-\varepsilon \mathbf{1},+\infty)$ cannot hold simultaneously.
Hence, $q=0$and so $pc\leq p(c_0-\varepsilon
\mathbf{1})$. The latter implies $(p(c_0-\varepsilon
\mathbf{1})-pc)_-=0$, that is $pc_0\leq p(c_0-\varepsilon
\mathbf{1})$, since $c_0=\bigvee \Lambda_-$. Hence, $p=0$.
However, this contradicts with the choise of $p>0$.
\end{proof}

\begin{lemma}
\label{l15}
The inequality $pe^a_z(-\infty,c_0+\varepsilon \mathbf{1})\succeq pe^a_z(c_0+\varepsilon \mathbf{1},+\infty)$ holds for all $\varepsilon>0$ and $0\neq p\in P(Z(\mathcal{M}))$.
\end{lemma}
\begin{proof}
Suppose a contrary. Then by Theorem \ref{t2} there exists a projection $0<p\in P(Z(\mathcal{M}))$, such that
for every $q\in P(Z(\mathcal{M})),\ 0<q\leq p$ the inequality
$qe^a_z(-\infty,c_0+\varepsilon \mathbf{1})\prec
qe^a_z(c_0+\varepsilon \mathbf{1},+\infty)$ holds.  Let $c\in
\Lambda_-$. Put $c_1=c(\mathbf{1}-p)+ (c_0+\varepsilon
\mathbf{1})p$. Then $c_1\in\Lambda_-.$ Hence, there exists a projection $0<r\in P(Z(\mathcal{M}))$ such that $re^a_z(-\infty,c_1)\succeq
re^a_z(c_1,+\infty).$ Thus, $r(\mathbf{1}-p)e^a_z(-\infty,c)\succeq r(\mathbf{1}-p)e^a_z(c,+\infty)$ and, therefore,  $r(\mathbf{1}-p)=0$. Also, $rpe^a_z(-\infty,c_0+\varepsilon \mathbf{1})\succeq rpe^a_z(c_0+\varepsilon \mathbf{1},+\infty)$ and, therefore, $rp=0$. It follows that $r=r(\mathbf{1}-p)+rp=0$. Hence, $c_1\in\Lambda_-$. However, $c_1\nleqslant c_0$. This contradicts to the definition of $c_0$.
\end{proof}

Let us consider the set

\begin{align*}
P_0:=& \{p\in P(Z(\mathcal{M})):\ \exists
c_p\in Z_h(LS(\mathcal{M})),\ q_p,r_p\in P(\mathcal{M}):\\ & q_p,r_p\leq
pe^a_z\{c_p\},\ q_pr_p=0,\ pe^a_z(-\infty,c_p)+q_p \sim pe^a_z(c_p,+\infty)+r_p\}
\end{align*}
and set
$$p_0:=\bigvee P_0.$$

\begin{lemma}
\label{l18} $P_0$ is a $\vee$-ideal in the Boolean algebra
$P(Z_h(\mathcal{M}))$ and $p_0\in P_0$, that is
$P_0=p_0P(Z_h(\mathcal{M}))$.
\end{lemma}
\begin{proof}

Let $p\in P_0,\ p>e\in P(Z_h(\mathcal{M}))$ and let $c_p,q_p,r_p$ are exactly the same like in the definition of $P_0$.
Then by letting  $q_e:=q_pe,r_e:=r_pe,c_e:=c_p$  one can see that these projections satisfy the condition in the definition of $P_0$.
Indeed, we have $q_er_e=eq_pr_pe=0,\ q_e,p_e\leq epe^a_z\{c_p\}=ee^a_z\{c_e\},\ ee^a_z(-\infty,c_e)+q_e=e(pe^a_z(-\infty,c_p)+q_p) \sim e(pe^a_z(c_p,+\infty)+r_p\}) =ee^a_z(c_e,+\infty)+r_e\}$. Hence, $e\in P_0$.

Let $e_1,e_2\in P_0$ и $e_1e_2=0$. Due to the definition of $P_0$, we know of the existence of $q_{e_i},r_{e_i}\in P(\mathcal{M})$, $c_{e_i}\in Z_h(LS(\mathcal{M}))$
such that $q_{e_i},r_{e_i}\leq e_i e^a_z\{c_{e_i}\},\ q_{e_i}r_{e_i}=0,\ e_i e^a_z(-\infty,c_{e_i})+q_{e_i} \sim e_i e^a_z(c_{e_i},+\infty)+r_{e_i}$
$i=1,2$. Let $c_{e_1+e_2}=c_{e_1}e_1+c_{e_2}e_2$. Then $c_{e_1+e_2}e_i = c_ie_i,\ i=1,2$.
So, we have $q_{e_1}+q_{e_2},r_{e_1}+r_{e_2}\leq (e_1+e_2)e^a_z\{c_{e_1+e_2}\},\ (q_{e_1}+q_{e_2})(r_{e_1}+r_{e_2})=0,\
(e_1+e_2) e^a_z(-\infty,c_{e_1+e_2})+(q_{e_1}+q_{e_2}) \sim (e_1+e_2) e^a_z(c_{e_1+e_2},+\infty)+(r_{e_1}+r_{e_2})$. So, $e_1+e_2\in P_0$.

Let now $e_1,e_2\in P_0$. Then $e_1\vee e_2=e_1+e_2(\mathbf{1}-e_1) \in P_0$, since $e_1[e_2(\mathbf{1}-e_1)]=0$
and $e_2(\mathbf{1}-e_1)\in P_0$. Hence, $P_0$ is a $\vee$-ideal in the Boolean algebra $P(Z_h(\mathcal{M}))$.

It follows from Lemma \ref{l16} that there exists a family $\{e_i\}_{i\in I}$ of pairwise disjoint projections from $P_0$
such that $\bigvee_{i\in I}e_i=p_0$. Since, $e_i\in P_0$, there are elements $c_{e_i},q_{e_i},r_{e_i}$. It follows from Lemma \ref{l17} that
there exists $c_{p_0}\in Z_h(LS(\mathcal{M}))$ such that
$c_{p_0}\mathbf{s}(c_{e_i})=c_{e_i}$ for every $i\in I$. Let
$q_{p_0}=\bigvee_{i\in I}q_{e_i},r_{p_0}=\bigvee_{i\in I}r_{e_i}$.
From the definition of projections $q_{p_0},r_{p_0}$ above, we have
$$q_{p_0},r_{p_0}\leq \bigvee_{i\in I}e_i e^a_z\{c_{e_i}\}=
p_0e^a_z\{c_{p_0}\}.
$$
Furthermore, since $e_ie_j=0$ for all $i\neq j$, we also have $q_{e_i}q_{e_j}=q_{e_i}r_{e_j}=0$ and therefore
$$
 q_{p_0}r_{p_0}=\bigvee_{i\in I}q_{e_i}
\bigvee_{i\in I}r_{e_i}=\bigvee_{i\in I}q_{e_i}r_{e_i}=0.
$$
Thus, we have
\begin{align*}
p_0e^a_z(-\infty,c_{p_0})+q_{p_0}&=\bigvee_{i\in
I}[e_ie^a_z(-\infty,c_{e_i})+q_{e_i}] \\ & \sim \bigvee_{i\in
I}[e_ie^a_z(c_{e_i},+\infty)+r_{e_i}]=p_0e^a_z(c_{p_0},+\infty)+r_{p_0},
\end{align*}
in other words, $p_0\in P_0$.
\end{proof}

\begin{lemma}
\label{l3} Suppose $p_0=\mathbf{1}$. Then there exists an element $u\in
U(\mathcal{M})$ such that
$|[a,u]|=u^*|a-c|u+|a-c|,\
u^2=\mathbf{1}$, where $c=c_{p_0}\in Z_h(LS(\mathcal{M}))$ is from the definition of the set $P_0$ for the element $p_0$.
\end{lemma}
\begin{proof}
Set $p=q_{p_0},q=r_{p_0}$ and $r:=\mathbf{1}-(e^a_z(-\infty,c)+p+
e^a_z(c,+\infty)+q)$. Then $p,q,r\leq e^a_z\{c\}$ and so
$ap=cp,\ aq=cq,\ ar=cr$. We claim that there exists a self-adjoint unitary element $u$
such that $u(e^a_z(-\infty,c)+p)=(e^a_z(c,+\infty)+q)u,\ ur=r$. Indeed, since $e^a_z(-\infty,c)+p\sim
e^a_z(c,+\infty)+q$, there exists a partial isometry  $v$ such that
$v^*v=e^a_z(-\infty,c)+p,\ vv^*=e^a_z(c,+\infty)+q$. Set $u:=v+v^*+r$. We have
$u^*u=e^a_z(-\infty,c)+p+e^a_z(c,+\infty)+q+r=\mathbf{1}$,
$uu^*=e^a_z(c,+\infty)+q+e^a_z(-\infty,c)+p+r=\mathbf{1}$,
$u^*=v^*+v+r=u$. This establishes the claim. It now remains to verify that \eqref{first main} holds.

To this end, first of all observe that the operators $a$ and $u^*au$ commute with the
projections $e^a_z(-\infty,c)+p$, $e^a_z(c,+\infty)+q$ and
$r$. This observation guarantees that
\begin{align*}
(u^*au-a)(e^a_z(-\infty,c)+p)&=|u^*au-a|(e^a_z(-\infty,c)+p),\\
(a-u^*au)(e^a_z(c,+\infty)+q)&=|u^*au-a|(e^a_z(c,+\infty)+q)
\end{align*}
and so
\begin{align*}
|u^*au-a|(e^a_z(-\infty,c)+p)&=
u^*a(e^a_z(c,+\infty)+q)u-a(e^a_z(-\infty,c)+p)\\
&=u^*a(e^a_z(c,+\infty)+q)u-c
u^*(e^a_z(c,+\infty)+q)u\\
&+
c(e^a_z(-\infty,c)+p)-a(e^a_z(-\infty,c)+p)\\
&=u^*|a(e^a_z(c,+\infty)+q)-c(e^a_z(c,+\infty)+q)|u\\
&+|c(e^a_z(-\infty,c)+p)-a(e^a_z(-\infty,c)+p)|\\
&=u^*|a-c|u(e^a_z(-\infty,c)+p)+
|a-c|(e^a_z(-\infty,c)+p).
\end{align*}
Similarly,
\begin{align*}
|u^*au-a|(e^a_z(c,+\infty)+q)&=
-u^*a(e^a_z(-\infty,c)+p)u+a(e^a_z(c,+\infty)+q)\\
&=-(u^*a(e^a_z(-\infty,c)+p)u-c
u^*(e^a_z(-\infty,c)+p)u)\\
&- c (e^a_z(c,+\infty)+q)
+a(e^a_z(c,+\infty)+q)\\
&=u^*|a-c|u(e^a_z(c,+\infty)+q)+
|a-c|(e^a_z(c,+\infty)+q).
\end{align*}
Finally, $(u^*au-a)r=c r - c r=0$, that is,
$|u^*au-a|r=0$.
We now obtain \eqref{first main} as follows

\begin{align*}|u^*au-a|&=|u^*au-a|[(e^a_z(-\infty,c)+p)+(e^a_z(c,+\infty)+q)+r]\\ & =
|u^*au-a|(e^a_z(-\infty,c)+p)+|u^*au-a|(e^a_z(c,+\infty)+q)+|u^*au-a|r\\ & =
(u^*|a-c|u+|a-c|)[(e^a_z(-\infty,c)+p)+(e^a_z(c,+\infty)+q)+r]\\ & =
u^*|a-c|u+|a-c|.
\end{align*}
\end{proof}

It follows from Lemma \ref{l3} that if $p_0=\mathbf{1}$, then
the part (i) of Theorem \ref{main} holds.

\begin{lemma}
\label{l31} Let $p\in P(Z(\mathcal{M}))$ be a finite projection.
Then $p\leq p_0$.
\end{lemma}
\begin{proof}
The case  $p=0$ is trivial. So, we assume $p>0$.
The algebra $\mathcal{M} p$ has a faithful normal centervalued
trace $\tau$ such that $\tau(p)=p$ \cite{KR2}. Then
$\tau(pe^a_z(-\infty,c))<\tau(pe^a_z(c,+\infty))$ for every
$c\in\Lambda_-$. Since $\tau$ is normal, using Lemma \ref{l11}
(vi),(vii) we have $\tau(pe^a_z(-\infty,c_0))=\tau(pe^a_z(-\infty,\bigvee \Lambda_-))=\tau(\bigvee pe^a_z(-\infty,\Lambda_-))=\bigvee \tau(pe^a_z(-\infty,\Lambda_-))\leq
\bigvee\tau(pe^a_z(\Lambda_-,+\infty))=\tau(\bigvee pe^a_z(\Lambda_-,+\infty))=\tau(pe^a_z(\bigvee\Lambda_-,+\infty))=
\tau(pe^a_z(c_0,+\infty))\leq\tau(pe^a_z[c_0,+\infty))$. Hence,
$\tau(pe^a_z(-\infty,c_0))\leq \tau(p)/2=p/2$. The same arguments with help of Lemma \ref{l15}
yield $\tau(pe^a_z(c_0+\varepsilon\mathbf{1},+\infty))\leq
\tau(p)/2=p/2$ for every $\varepsilon>0$. Since
$e^a_z(c_0,+\infty))=\bigvee_{\varepsilon>0}
e^a_z(c_0+\varepsilon\mathbf{1},+\infty)$, we have
$$
\tau(pe^a_z(c_0,+\infty))\leq p/2.
$$
By Theorem \ref{t2}, there exist projections $p_1,p_2\in P(Z(\mathcal{M}))$, such that
$$p_1p_2=0,\ p_1+p_2=p,\ p_1 e^a_z(-\infty,c_0)\preceq p_1
e^a_z(c_0,+\infty),\ p_2 e^a_z(-\infty,c_0)\succeq p_2
e^a_z(c_0,+\infty).
$$
Applying then the inequality above, we have
\begin{align*}
p_1\ge 2\tau(p_1e^a_z(c_0,+\infty))=&\tau(p_1e^a_z(c_0,+\infty))+\tau(p_1e^a_z(-\infty,c_0))+
(\tau(p_1e^a_z(c_0,+\infty))\\ &-\tau(p_1e^a_z(-\infty,c_0))),
\end{align*}
and, immediately,
\begin{align*}\tau(p_1e^a_z(c_0,+\infty))-\tau(p_1e^a_z(-\infty,c_0))&\leq
\tau(p_1(\mathbf{1}-e^a_z(-\infty,c_0)-e^a_z(c_0,+\infty))) \\ &=
\tau(p_1 e^a_z\{c_0\}).
\end{align*}
Hence,
$$p_1e^a_z(-\infty,c_0)\preceq p_1e^a_z(c_0,+\infty)\preceq
p_1e^a_z(-\infty,c_0)+p_1 e^a_z\{c_0\}.
$$ Hence, it follows from Lemma
\ref{l00} that there exists a projection
$P(\mathcal{M})\ni p_{11}\leq p_1 e^a_z\{c_0\}$ such that
$p_{11}+p_1e^a_z(-\infty,c_0)\sim p_1e^a_z(c_0,+\infty)$.

Analogous arguments show that there exists a projection
$P(\mathcal{M})\ni p_{21}\leq p_2 e^a_z\{c_0\}$, such that
$p_2e^a_z(-\infty,c_0)\sim p_2e^a_z(c_0,+\infty) + p_{21}$.

Since, $p_1+p_2=p$ и $p_1p_2=0$ we have $p_{11}+pe^a_z(-\infty,c_0)\sim pe^a_z(c_0,+\infty) +
p_{21}$. Hence, $p\in P_0$.
\end{proof}

In the proof of the following lemma we shall frequently use a well known fact that $\mathbf{c}(q)z=\mathbf{c}(qz),\ \forall q\in P(\mathcal{M}),\ \ \forall z\in P(Z(\mathcal{M}))$ \cite[Prop. 5.5.3]{KR}.

\begin{lemma}
\label{l32} Let $p\in P(Z(\mathcal{M}))$ and let $\mathcal{M}p$ be
a $\sigma$-finite purely infinite $W^*$-subalgebra in
$\mathcal{M}$. Then $p\leq p_0$.
\end{lemma}
\begin{proof}
Suppose $p\nleqslant p_0$. Then, without loss of generality we can assume $p p_0=0$.
Indeed, since $p=pp_0+p(\mathbf{1}-p_0)$, we see that $p\nleqslant
p_0\Leftrightarrow p(\mathbf{1}-p_0)>0$. Assuming
$p(\mathbf{1}-p_0)>0$, we observe that the algebra
$\mathcal{M}p(\mathbf{1}-p_0)$ is a $\sigma$-finite purely
infinite $W^*$-subalgebra in $\mathcal{M}$. Therefore, if we show
that $p(\mathbf{1}-p_0)\in P_0$, that is $p(\mathbf{1}-p_0)\leq
p_0$, then we would also have that $p\in P_0$.

Applying Theorem \ref{t2}
to $e^a_z(-\infty,c_0)$  and $e^a_z(c_0,+\infty)$ and setting
$p_1=z_-,\ p_2=z_+$ and multiplying $z_-,z_+$ by $p$, we obtain $p_1p_2=0$, $p_1+p_2\leq
p$ and so
$$p_1
e^a_z(-\infty,c_0)\prec\prec p_1 e^a_z(c_0,+\infty),\ p_2
e^a_z(-\infty,c_0)\succ\succ p_2 e^a_z(c_0,+\infty),
$$
and
$$(p-p_1-p_2)e^a_z(-\infty,c_0)\sim (p-p_1-p_2)e^a_z(c_0,+\infty).$$
Then $p-p_1-p_2\in P_0$, and so $p-p_1-p_2\leq p_0$. Since,
$pp_0=0$ and $p_1+p_2\leq p$, we have $p_1+p_2=p$.
Next, we  define the following three
projections $q_-,q_0,q_+\in P(Z(\mathcal{M}))$, by setting
$$q_-:=\mathbf{c}(e^a_z(-\infty,c_0))(\mathbf{1}-\mathbf{c}(e^a_z(c_0,+\infty))),$$
$$q_+:=\mathbf{c}(e^a_z(c_0,+\infty))(\mathbf{1}-\mathbf{c}(e^a_z(-\infty,c_0)))),\
q_0:=\mathbf{c}(e^a_z(-\infty,c_0))\mathbf{c}(e^a_z(c_0,+\infty)).$$
It follows from the definition that
$$q_- q_0 = q_0 q_+=q_+ q_-=0,\ q_-
+q_0+q_+=\mathbf{c}(e^a_z(-\infty,c_0))\vee
\mathbf{c}(e^a_z(c_0,+\infty)).$$
Next, we have for $i=1,2$
$$
\mathbf{c}(q_0p_ie^a_z(-\infty,c_0))=q_0p_i\mathbf{c}(e^a_z(-\infty,c_0))=$$
$$=q_0p_i=q_0p_i\mathbf{c}(e^a_z(c_0,+\infty))=\mathbf{c}(q_0p_ie^a_z(c_0,+\infty)).$$
Hence, $$q_0p_ie^a_z(-\infty,c_0)\sim
q_0p_ie^a_z(c_0,+\infty)$$ and $$q_0pe^a_z(-\infty,c_0)\sim
q_0pe^a_z(c_0,+\infty)$$ (see \cite{Sak} Proposition 2.2.14).
The preceding estimates imply that
$q_0p\in P_0$ (indeed, it is sufficient to set $c_{q_0 p}=c_0,\
r_{q_0 p}=q_{q_0 p}=0$), and hence $q_0p\leq p_0$.
Since, $q_0p\leq p$ and $pp_0=0$ we have $q_0p=0$. It is clear that
$p_1q_-=0$ (or $p_1q_-e^a_z(-\infty,c_0)\prec
p_1q_-e^a_z(c_0,+\infty)=0$). Analogously, $p_2q_+=0$.

From Lemma \ref{l15}, we have $pe^a_z(-\infty,c_0+\varepsilon
\mathbf{1})\succ\succ pe^a_z(c_0+\varepsilon \mathbf{1},+\infty)$
for every $\varepsilon>0$ (if $P(Z(\mathcal{M}))\ni q\leq p$ and
$qe^a_z(-\infty,c_0+\varepsilon \mathbf{1})\sim
qe^a_z(c_0+\varepsilon \mathbf{1},+\infty)$ then $q\in P_0$.
Hence, $q=0$). Applying Lemma \ref{l111} to the element $ap$ in
the algebra $\mathcal{M}p$, we obtain
$\mathbf{c}(pe^a_z(-\infty,c_0+\varepsilon \mathbf{1}))=p$. Then
$pe^a_z(c_0+\varepsilon \mathbf{1},+\infty)=0$ for every
$\varepsilon>0$. Setting $r:=pe^a_z(c_0+\varepsilon
\mathbf{1},+\infty)>0$ and observing that $pe^a_z(-\infty,c_0+\varepsilon
\mathbf{1})\succ\succ pe^a_z(c_0+\varepsilon \mathbf{1},+\infty)$ and $\mathbf{c}(r)p=\mathbf{c}(pe^a_z(c_0+\varepsilon
\mathbf{1},+\infty))$, we have
\begin{equation}\label{A}
\mathbf{c}(r)pe^a_z(-\infty,c_0+\varepsilon \mathbf{1})
\succ \mathbf{c}(r)pe^a_z(c_0+\varepsilon \mathbf{1},+\infty).
\end{equation}
Next, we observe that
\begin{align*}\mathbf{c}(\mathbf{c}(r)pe^a_z(-\infty,c_0+\varepsilon
\mathbf{1}))&=\mathbf{c}(r)\mathbf{c}(pe^a_z(-\infty,c_0+\varepsilon
\mathbf{1}))\\ & =\mathbf{c}(r)p=\mathbf{c}(rp)=\mathbf{c}(r).
\end{align*}
Thus, we have
$$\mathbf{c}(\mathbf{c}(r)pe^a_z(-\infty,c_0+\varepsilon
\mathbf{1}))=\mathbf{c}(r)=\mathbf{c}(\mathbf{c}(r)pe^a_z(c_0+\varepsilon
\mathbf{1},+\infty))$$
which however is a contradiction with \eqref{A} in view of \cite[Proposition 2.2.14]{Sak}.

Hence,
$pe^a_z(c_0,+\infty)=\bigvee_{\varepsilon>0}pe^a_z(c_0+\varepsilon
\mathbf{1},+\infty)=0$. Since $p_1\leq p$, we have $p_1e^a_z(c_0,+\infty)=0$.
Using the assumption
$p_1e^a_z(-\infty,c_0)\prec\prec p_1e^a_z(c_0,+\infty)$, we obtain
$p_1=0$.

Next, for every $c\in \Lambda_-$ we have by the definition of $\Lambda_-$ that either $p_2q_-e^a_z(-\infty,c)=0$, or else
$p_2q_-e^a_z(-\infty,c)\prec\prec p_2q_-e^a_z(c,+\infty)$.
In the first case, taking the supremum and appealing to \cite[Prop. 5.5.3]{KR}, we have
$$p_2q_-=p_2q_-\mathbf{c}(e^a_z(-\infty,c_0))=\mathbf{c}(p_2q_-e^a_z(-\infty,c_0))=0.$$
In the second case, assuming $p_2q_-\neq 0$,
we obtain that $\mathcal{M} p_2q_-$ is a purely infinite
$\sigma$-finite $W^*$-algebra. Hence, using the same argument as in the previous paragraph,
we have $p_2q_-e^a_z(-\infty,c_0)=0$. So, we have obtained a contradiction with
the choice of $p_2$. Indeed, we assumed that $p_2
e^a_z(-\infty,c_0)\succ\succ p_2 e^a_z(c_0,+\infty)$ and defined $q_-$ so that  $p_2q_-\leq p_2$. Hence,
$p_2q_-=0$. Since we have already shown that $p=p_1+p_2$, we can now write $$p(\mathbf{c}(e^a_z(-\infty,c_0))\vee
\mathbf{c}(e^a_z(c_0,+\infty)))=(p_1+p_2)(q_-+q_0+q_+)=0.$$

Let $q=pe^a_z\{c_0\}$. Then $aq=c_0q$ and
$qe^a_z(-\infty,c_0)=qe^{aq}_z(-\infty,qc_0)=\mathbf{s}((c_0q-aq)_+)=0$
and, in the same way, $qe^a_z(c_0,+\infty)=0$. Hence, $q\in P_0$ and
$q\leq p_0$. Since $q\leq p$ and
$pp_0=0$, we infer $q=0$. Again appealing to \cite[Prop. 5.5.3]{KR}, we obtain $p\mathbf{c}(e^a_z\{c_0\})=0$.  Hence,
$p=p(\mathbf{c}(e^a_z(-\infty,c_0))\vee
\mathbf{c}(e^a_z(c_0,+\infty)))\vee \mathbf{c}(e^a_z\{c_0\}))=0$.
Then $\mathcal{M}p=\{0\}$, which contradicts with the lemma's assumption.
Hence, $p\leq p_0$.
\end{proof}

Let $q\in P(Z(\mathcal{M}))$ and suppose
$$qe^a_z(-\infty,c_0]=qe^a_z(-\infty,c_0)+qe^a_z\{c_0\}\sim qe^a_z(c_0,+\infty).
$$ Then setting $c_q:=c_0,\ q_q:=qe^a_z\{c_0\},\ r_q:=0$ we see that $q\in P_0$ and therefore
$q\leq p_0$.
So, there exist projections $p_+$ и $p_-\in
P(Z(\mathcal{M}))$ (which may be null projections), such that
$$p_-+p_+=\mathbf{1}-p_0,\ p_-p_+=0,
$$
and
$$
p_+e^a_z(-\infty,c_0]\succ\succ p_+e^a_z(c_0,+\infty),\
p_-e^a_z(-\infty,c_0]\prec\prec p_-e^a_z(c_0,+\infty).
$$ By Lemma
\ref{l31}, the following implications hold \\$p_-\neq 0$ (resp.
$p_+\neq 0$) $\Rightarrow$ $p_-\mathcal{M}$ (resp.
$p_+\mathcal{M}$) is a properly infinite $W^*$-algebra.

\begin{lemma}
\label{l33} For every $q\in P(Z(\mathcal{M})),\ 0<q\leq p_+$, we have $qe^a_z(-\infty,c_0)\succ qe^a_z[c_0,+\infty)$.
\end{lemma}
\begin{proof}
Firstly, we show that $p_+e^a_z(-\infty,c_0)\succeq
p_+e^a_z(c_0,+\infty)$. If this fails then there exists an element $0\neq q\in
P(Z(\mathcal{M}p_+))$ such that
$$qe^a_z(-\infty,c_0) \prec
qe^a_z(c_0,+\infty).$$
On the other hand, we have
$$qe^a_z(-\infty,c_0]\succ\succ qe^a_z(c_0,+\infty).$$
Indeed, to see the preceding estimate, let  $0<r\leq q$. Then we
have $r\leq p_+=
p_+e^a_z(-\infty,c_0]\vee
p_+e^a_z(c_0,+\infty)$, and it
follows from the definition of the symbol "$\succ\succ$" that
$re^a_z(-\infty,c_0]\succ re^a_z(c_0,+\infty)$.
By Lemma \ref{l00}, we have that $qe^a_z(c_0,+\infty)\sim
qe^a_z(-\infty,c_0)+r$, where $r\in P(Z(\mathcal{M})q)$ and
$r<qe^a_z\{c_0\}$. In this case, setting $c_q:=c_0,\ q_q:=r,\
r_q:=0$ we obtain $q\in P_0$ and so $q\leq p_0$. Thus, $q\leq
p_0p_+=0$. This contradiction shows that
\begin{equation}\label{B}
p_+e^a_z(-\infty,c_0)\succeq
p_+e^a_z(c_0,+\infty).
\end{equation}

Let us now consider projections $p_+e^a_z(-\infty,c_0)$ and
$p_+e^a_z[c_0,+\infty)$ in the algebra $p_+\mathcal{M}$. In the
notation of Theorem \ref{t2} applied to the algebra
$p_+\mathcal{M}$, we intend to prove that $z_-=z_0=0$ and so
$z_+=p_+$.  Suppose that there exists $0<r\leq p_+$, $r\in
P(Z(\mathcal{M}))$ such that $re^a_z(-\infty,c_0)\sim
re^a_z[c_0,+\infty)=re^a_z\{c_0\}+re^a_z(c_0,+\infty)$. Then
$r\leq p_0$ and therefore $r\leq p_+p_0=0$. This shows that
$z_0=0$. Next, we shall show the equality $z_-=0.$ Supposing that
$z_->0$, we have $z_-e^a_z(-\infty,c_0)\prec
z_-e^a_z[c_0,+\infty)$.  Then $z_-e^a_z[c_0,+\infty)=p_1+p_2$, where
$p_1,p_2\in P(\mathcal{M}),\ p_1p_2=0,\ p_1\sim
z_-e^a_z(-\infty,c_0)$. From $p_+e^a_z(-\infty,c_0)\succeq
p_+e^a_z(c_0,+\infty)$ it follows that
$z_-e^a_z(-\infty,c_0)\succeq z_-e^a_z(c_0,+\infty)$. Then, by
Lemma \ref{l00}, we know that there exists some $q\in
P(\mathcal{M}z_-)$, such that $z_-e^a_z(-\infty,c_0)\sim
z_-e^a_z(c_0,+\infty)+q$ and $q<e^a_z\{c_0\}z_-$. Hence, $z_-\in
P_0$, i.e. $z_-\leq p_0$. Therefore, $z_-p_+=0$ and $z_-=0$.
This contradiction completes the proof.
\end{proof}

\begin{lemma}
\label{l33_1} $qe^a_z(c_0,c_0+\varepsilon \mathbf{1}]\succ
qe^a_z(-\infty,c_0]+qe^a_z(c_0+\varepsilon \mathbf{1},+\infty)$
for every $\varepsilon>0$ and every $q\in P(Z(\mathcal{M})),\
0<q\leq p_-$.
\end{lemma}
\begin{proof}
Suppose, there exists $0<p\in P(Z(\mathcal{M}p_-))$ such that
$pe^a_z(c_0,c_0+\varepsilon \mathbf{1}]\preceq
pe^a_z(-\infty,c_0]+pe^a_z(c_0+\varepsilon \mathbf{1},+\infty)$.
By the assumption the projection $p_-\neq 0$ and hence it is properly infinite (see the implication preceding
Lemma \ref{l33}) and since $p\leq p_-$, we conclude
that $p$ is also properly infinite. Due to Lemma \ref{l0} (iii),
we have $pe^a_z(-\infty,c_0]+pe^a_z(c_0+\varepsilon
\mathbf{1},+\infty)\sim p$. However,
$pe^a_z(-\infty,c_0]\prec\prec pe^a_z(c_0,+\infty)$. Hence,
$pe^a_z(-\infty,c_0]\prec\prec p$, and by Lemma \ref{l0} (ii)
 we have $pe^a_z(c_0+\varepsilon
\mathbf{1},+\infty)\sim p$ (indeed, otherwise, we would have had
$pe^a_z(-\infty,c_0]+pe^a_z(c_0+\varepsilon
\mathbf{1},+\infty)\prec p$). Next, it follows from Lemma
\ref{l15_0} that $pe^a_z(-\infty,c_0+\varepsilon
\mathbf{1})\succeq pe^a_z(c_0+\varepsilon \mathbf{1},+\infty)$,
that is $pe^a_z(-\infty,c_0+\varepsilon \mathbf{1})\sim p \sim
pe^a_z(c_0+\varepsilon \mathbf{1},+\infty)$. Then $p\leq p_0$.
Which is a contradiction with the assumption $0<p\leq p_-$.
\end{proof}

In the case when $p_+=\mathbf{1}$, it follows from Lemma \ref{l33} that $c_0\in\Lambda_+$ and therefore in this case
$\bigvee\Lambda_-=c_0=\bigwedge\Lambda_+$, where the first equality is simply the definition of $c_0$.
Let us explain the second equality. By Lemma \ref{l14}, we have $c_1<c_2$ for every $c_1\in\Lambda_-,\ c_2\in\Lambda_+$. Therefore, $c_0\leq c_2$ for any $c_2\in\Lambda_+$.
However, $c_0\in\Lambda_+$, and hence $c_0=\bigwedge\Lambda_+$.

In the case when $p_-=\mathbf{1}$, appealing to the definitions of $p_-$ and $\Lambda_-$ we have $c_0\in\Lambda_-$.
Moreover,
in this case we have by Lemma \ref{l33_1} that
$e^a_z(c_0+2\varepsilon\mathbf{1},+\infty)\leq
e^a_z(c_0+\varepsilon\mathbf{1},+\infty)\prec\prec
e^a_z(c_0,c_0+\varepsilon\mathbf{1}]\leq
e^a_z(c_0,c_0+2\varepsilon\mathbf{1})\leq
e^a_z(-\infty,c_0+2\varepsilon\mathbf{1})$, that is
$c_0+2\varepsilon\mathbf{1}\in \Lambda_+$ for every
$\varepsilon>0$. Then it follows from Lemma \ref{l11} (v) that
in the case $p_-=\mathbf{1}$, we have
$\bigvee\Lambda_-=c_0=\bigwedge\Lambda_+$ as well. Indeed,
$c_0=\bigwedge_{\varepsilon>0}(c_0+\varepsilon\mathbf{1})\geq\bigwedge\Lambda_+\geq c_0$.

Hence, the cases
$p_+=\mathbf{1}$ and $p_-=\mathbf{1}$ are symmetric, that is the second case may be obtained from the first case using
the substitution $a\rightarrow -a$ and $c_0\rightarrow -c_0$.
Hence, in the sequel we will consider only the case
$p_-=\mathbf{1}$. In this case, the algebra $\mathcal{M}$ is properly infinite, since the projection $p_-$ is properly infinite.
It follows from Lemma \ref{l33_1} that in this case
\begin{equation}\label{last condition}qe^a_z(c_0,c_0+\varepsilon \mathbf{1}]\succ
qe^a_z(-\infty,c_0]+qe^a_z(c_0+\varepsilon \mathbf{1},+\infty)\end{equation}
for every $\varepsilon>0$ and  every  $0<q\in P(Z(\mathcal{M}))$.

The following lemma extends the result of \cite[Lemma 4]{BS1}.

\begin{lemma}
\label{l1} Let $p,q\in P(\mathcal{M})$, $p\succeq q$ and let
one of the following holds:

(i). $q$ is a finite projection and there exists a non-decreasing
sequence $\{p_n\}$ of finite projections in $\mathcal{M}$ such
that $p_n \uparrow p$ and $ap_n=p_n a$ for every $n\in \mathbb{N}$;

(ii). $q$ is a properly infinite projection and $ap=pa\in \mathcal{M}$.

Then there exist a projection $q_1\leq p$ in $\mathcal{M}$ such that
$q_1\sim q$ and $aq_1=q_1a$.
\end{lemma}
\begin{proof}
Assume (i) holds.

Set $\mathcal{A}_1:=\{b\in\mathcal{M}:\ ba=ab\}$ (that is
$\mathcal{A}_1$ is the commutant of the family of all spectral projections of $a$).
Since, $p_n \uparrow p$ and $ap_n=p_n a$ for every $n\in
\mathbb{N}$, we have $ap=pa$. Let $\mathcal{A}:=p\mathcal{A}_1$.
Then $\mathcal{A}$ is $W^*$-subalgebra in $\mathcal{M}$ with identity
$p$.

First of all, let us show that every atom of algebra $\mathcal{A}$ (if
it exists) is an atom of the algebra $\mathcal{M}$. Let
$e$ be an atom of the algebra $\mathcal{A}$, $0\leq f<e,\ f\in
P(\mathcal{M})$. So, if $g$ is a spectral projection of $a$ then
$gf=g(p(ef))=((gp)e)f\in \{0,e\}f=\{0,ef\}=\{0,f\}\subset
\mathcal{M}_h$. In particular, $gf=fg$, that is $f\in\mathcal{A}$.
Hence, $f=0$. Consequently, $e$ is an atom of $\mathcal{M}$.

Let $z_n=\bigvee\{z\in P(Z(\mathcal{M})):\ zq\preceq zp_n\}$.
It is easy to see, that $z_n\uparrow\mathbf{1}$ (otherwise
$(\mathbf{1}-\bigvee_n z_n)q \succ (\mathbf{1}-\bigvee_n z_n)p$).
We will assume $z_nq>0$ for every $n\in \mathbb{N}$.

Let us construct a non-decreasing sequence of projections
$\{f_n\}$ in $\mathcal{A}$ such $p_nz_n\geq f_n\sim qz_n$. Let
$f_0=z_0=0$. Suppose that $f_0,f_1,...,f_{n-1}$ have been constructed.

The set $\{r\in P(\mathcal{A}):\ qz_n\preceq r,\ f_{n-1}\leq r
\leq p_n z_n\}$ is non-empty (since, it contains $p_n z_n$) and is contained
in the finite $W^*$-algebra $(p_n\vee q)\mathcal{A} (p_n\vee q)$. Let
$\tau$ be a center-valued trace on the algebra $(p_n\vee q)\mathcal{A}
(p_n\vee q)$. From the Zorn's lemma and the fact that $\tau$ is normal, we have
that this set has the minimal element $r_0$. Suppose
 $\tau(r_0)>\tau(qz_n)$. Using the Zorn's lemma and the fact that $\tau$ is normal
again, we obtain that the set $\{r\in P(\mathcal{A}):\
\tau(r)<\tau(qz_n),\ 0 \leq r\leq r_0\}$ has a maximal element
 $r_1$. Then
$0\leq\tau(r_1)<\tau(qz_n)<\tau(r_0)$.

We claim that $r_0-r_1$ is an atom in the algebra $\mathcal{A}$.
Suppose there exists $e\in P(\mathcal{A}),\ 0\leq
e<r_0-r_1$. Then there exists a central projection $z$ in the
algebra $\mathcal{M}$ such that $z(r_1+e)\preceq z(q z_n),\
(\mathbf{1}-z)(r_1+e)\succeq (\mathbf{1}-z)(q z_n)$. If $ze>0$
then $P(\mathcal{A})\ni z(r_1+e)+(\mathbf{1}-z)r_1>r_1,\
\tau(z(r_1+e)+(\mathbf{1}-z)r_1)\leq\tau(z(qz_n)+(\mathbf{1}-z)qz_n)=\tau(q z_n)$ and if
$(\mathbf{1}-z)e>0$ then $P(\mathcal{A})\ni
zr_0+(\mathbf{1}-z)(r_1+e)<r_0,\
\tau(zr_0+(\mathbf{1}-z)(r_1+e))\ge\tau(q z_n)$. The first assumption
contradicts the maximality of $r_1$, the second assumption
contradicts the minimality of $r_0$.
Hence, $ze=(\mathbf{1}-z)e=0$ that is $e=0$.
Hence, $r_0-r_1$ is an atom of algebra $\mathcal{A}$.

Hence, $r_0-r_1$ is an atom of algebra $\mathcal{M}$.

Since, $\tau(r_1)<\tau(qz_n)$, we have $r_1\prec qz_n$. Hence,
there exists a  projection $e_1\in (p_n\vee q)\mathcal{M} (p_n\vee
q)$ such that $r_1\sim e_1<qz_n$. Then
$\tau(qz_n-e_1)=\tau(qz_n)-\tau(r_1)<\tau(r_0-r_1)$, that is
$qz_n-e_1\prec r_0-r_1$. Since $r_0-r_1$ is an atom of $(p_n\vee
q)\mathcal{M} (p_n\vee q)$, we infer $qz_n-e_1=0$, which is contradicts to the choice of
$e_1$. Then $\tau(r_0)=\tau(qz_n)$, that is $r_0\sim qz_n$.
We set $f_n:=r_0$.

This completes the construction of the sequence $\{f_n\}_{n\geq 1}$.
Let $q_1:=\bigvee_{n=1}^\infty f_n$.

Since $f_n\sim qz_n,\ f_{n+1}\sim qz_{n+1}$ and all four
projections are finite we have $f_{n+1}-f_n\sim qz_{n+1}-qz_n$.
Indeed, applying the center-valued trace $\tau$ on the finite $W^*$-algebra
$(f_{n+1}\vee f_n\vee qz_{n+1}\vee qz_n)\mathcal{M}(f_{n+1}\vee
f_n\vee qz_{n+1}\vee qz_n)$ trivially yields
$\tau(f_{n+1}-f_n)=\tau(f_{n+1})-\tau(f_n)=\tau(qz_{n+1})-\tau(qz_n)=\tau(qz_{n+1}-qz_n)$. The latter implies immediately
$f_{n+1}-f_n\sim qz_{n+1}-qz_n$.

Hence, $q_1=f_1\vee\bigvee_{n=1}^\infty (f_{n+1}-f_n)\sim
qz_1\vee\bigvee_{n=1}^\infty (qz_{n+1}-qz_n)=q$.

Assume (ii) holds. By the assumption there exists a projection
$q_1^0\in \mathcal{M}$, such that $q_1^0\leq p$ and $q_1^0\sim q$.
We set
$$q_1^n:=\mathbf{l}(a^n q_1^0),\ \forall n>0,\
q_1:=\bigvee_{k=0}^\infty q_1^k.$$ We claim that $q_1\sim q$.
Indeed, since $q_1\ge q_1^0\sim q$, we have $q_1\succeq q$. On the
other hand, we have $q_1^n\sim \mathbf{r}(a^nq_1^0)\leq q_1^0\sim
q$, which implies $q_1^n \preceq q$ for all $n\geq 0$. Now, we
shall show that in fact $q_1 \preceq q$. Note that although  $q$
is a properly infinite projection we cannot simply refer to Lemma
\ref{l0}(i) since the sequence $\{q_1^k\}_{k\ge 0}$ does not
necessarily consist of pairwise orthogonal elements. However,
representing the projection $q_1$ as
\begin{align*}
q_1&=\bigvee_{k=0}^\infty q_1^k=\sum_{m=1}^\infty(\bigvee_{k=0}^m
q_1^k-\bigvee_{k=0}^{m-1} q_1^k)+q_1^0
=\sum_{m=1}^\infty(q_1^m\vee \bigvee_{k=0}^{m-1}
q_1^k-\bigvee_{k=0}^{m-1} q_1^k) +q_1^0.
\end{align*}
and noting that $q_1^0\sim q$ and
\begin{align*}
q_1^m\vee\bigvee_{k=0}^{m-1} q_1^k-\bigvee_{k=0}^{m-1} q_1^k\sim
q_1^m
-(q_1^m\wedge \bigvee_{k=0}^{m-1} q_1^k) &\leq q_1^m \preceq q
\end{align*}
we infer via Lemma \ref{l0}(i) that $q_1 \preceq q$.  This completes the proof of the claim.

Since $ap=pa$ and $q_1^0 \leq p$, we have $pa^n q_1^0 = a^n p
q_1^0=a^n q_1^0$, and so $q_1^n\leq p$ for all $n>0$. Hence,
$q_1\leq p$. It remains to show that $aq_1=q_1a$. The subspace
$q_1(H)$ coincides with the closure of linear span of the set $Q:=\{ a^n
q_1^0(H):\ n>0\}$. Indeed, setting $Q_k:=a^k q_1^0(H)$, we see that $q_1^k(H)=\overline{Q_k}$ and so
$\overline{Q}\supset q_1^k(H)$ for every $k>0$. Therefore,
$q_1(H)=\overline{\bigcup_{k=1}^\infty q_1^k(H)}\subset
\overline{Q}$. Conversely, $q_1^k(H)\subset q_1(H)$, for every $k>0$. Therefore, the closed linear span containing the set
$\overline{\bigcup_{k=1}^\infty q_1(H)}$ contains $\overline{Q}$ and is contained in
$q_1(H)$.

By the assumption the operator $ap$ is bounded,
and since $q_1\leq p$, the operator $aq_1$ is also bounded. Thus,
for every vector $\xi\in Q$, the vector $a\xi=aq_1\xi$ again
belongs to $Q$. Again appealing to the fact that $aq_1$ is
bounded, we infer $q_1 a q_1=a q_1$. From this we conclude that
$aq_1=q_1 a$.
\end{proof}

\begin{lemma}
\label{l20} Let $0<b\in Z(\mathcal{M})$,
$\mathbf{s}(b)=\mathbf{1}$; $e^a_z(\mathbf{0},\infty)$ be a
properly infinite projection and
$\mathbf{c}(e^a_z(\mathbf{0},\infty))=\mathbf{1}$. Let projection
$q\in P(\mathcal{M})$ be finite or properly infinite,
$\mathbf{c}(q)=\mathbf{1}$ and $q\prec\prec
e^a_z(\mathbf{0},\infty)$. Let $\mathbb{R}\ni \mu_n\downarrow
0$. For every $n\in \mathbb{N}$ we denote by $z_n$ such a projection
that $\mathbf{1}-z_n$ is the largest central projection,
for which $(\mathbf{1}-z_n)q\succeq (\mathbf{1}-z_n)e^a_z(\mu_n
b,+\infty)$ holds. We have $z_n\uparrow_n \mathbf{1}$ and for $$d:=[\mu_1 z_1
+\sum_{n=1}^\infty \mu_{n+1}(z_{n+1}-z_n)]b$$ the following relations hold:
$q\prec\prec e^a_z(d,+\infty)$, $0<d\leq \mu_1 b$ and
$\mathbf{s}(d)=\mathbf{1}$. Moreover, if all projections
$e^a_z(\mu_n b,+\infty)$, $n\geq1$ are finite then $e^a_z(d,+\infty)$ is a
finite projection as well.
\end{lemma}
\begin{proof}
Since, $e^a_z(\mu_{n+1} b,+\infty)\geq e^a_z(\mu_n b,+\infty)$
(by Lemma \ref{l11} (i)) we have $(\mathbf{1}-z_{n+1})q\succeq
(\mathbf{1}-z_{n+1})e^a_z(\mu_{n+1} b,+\infty)\geq
(\mathbf{1}-z_{n+1})e^a_z(\mu_n b,+\infty)$. Hence, $z_{n+1}\geq
z_n$ for every $n\in \mathbb{N}$. In addition, $e^a_z(\mu_n
b,+\infty)\uparrow_n e^a_z(\mathbf{0},+\infty)$ (by Lemma \ref{l11}
(vii)) and $e^a_z(\mathbf{0},+\infty)$ is properly infinite projection.
Hence, in the case when $q$ is finite projection, it follows from Lemma
 \ref{l0} (v) that $z_n\uparrow_n \mathbf{1}$. Let us consider the case when
$q$ is a properly infinite projection with
$\mathbf{c}(q)=\mathbf{1}$  and such that $q\prec\prec e^a_z(\mathbf{0},\infty)$. In this case, we apply Lemma \ref{l0}
(vi) with  $p=q,\ q=e^a_z(0,+\infty),\ q_n=e^a_z(\mu_n b,+\infty)$ and deduce  $\bigvee_{n=1}^\infty z_n\geq
\mathbf{c}(q)=\mathbf{1}$.

All other statements follow from the form of element $d$. Since, $z_1 d
=\mu_1 z_1 b$, $(z_{n+1}-z_n)d=\mu_{n+1}(z_{n+1}-z_n)b$ and $z_n
q\prec\prec z_n e^a_z(\mu_n b,+\infty)$ for every $n\in
\mathbb{N}$. Observe also that $\mathbf{s}(d)=\mathbf{s}(b)(z_1+\sum_{n=1}^\infty (z_{n+1}-z_n))=\mathbf{1}$.

Finally, let all projections $e^a_z(\mu_n b,+\infty)$, $n\geq1$ be finite. Since $dz_1=\mu_1 b,\ d(z_{n+1}-z_n)=\mu_{n+1}b(z_{n+1}-z_n)$, we have
$$e^a_z(d,+\infty)z_1=e^a_z(\mu_1 b,+\infty)z_1,
$$
$$ e^a_z(d,+\infty)(z_{n+1}-z_n)=e^a_z(\mu_{n+1} b,+\infty)(z_{n+1}-z_n)$$ for every $n\in\mathbb{N}$.
There projections standing on the right-hand sides are finite . Hence, $e^a_z(d,+\infty)$ is finite projection as a sum of the left-hand sides \cite[Lemma 6.3.6]{KR2}.
\end{proof}

In the proof of the following lemma, we shall use a following well known implication
$$p\prec\prec q \Rightarrow zp\prec\prec zq, \quad \forall
z\in P(Z(\mathcal{M})).$$
We supply here a straightforward argument for convenience of the reader. Let $z'\in z\in Z(\mathcal{M})$ be such that $0<z'\leq
\mathbf{c}(pz)\vee\mathbf{c}(qz)=z(\mathbf{c}(p)\vee\mathbf{c}(q))$.
Then $z'\leq\mathbf{c}(p)\vee\mathbf{c}(q)$ and therefore
$z'(zp)=z'p\prec z'q=z'(zq)$. This means $zp\prec\prec zq$.

\begin{lemma}
\label{l2} Let $\mathcal{M}$ be  properly infinite and
$$e^a_z(\mathbf{0},t \mathbf{1}]\succ\succ
e^a_z(-\infty,\mathbf{0}]+e^a_z(t\mathbf{1},+\infty)$$ for every
$t>0$. Then for every $\varepsilon>0$ there exists an element
$u_\varepsilon\in U(\mathcal{M})$ such that $|[a,u_\varepsilon]|\geq
(1-\varepsilon)|a|,\ u_\varepsilon^2=\mathbf{1}$.
\end{lemma}
\begin{proof}
Of course, we may assume $\varepsilon<1$.

Let us observe that
$\mathbf{c}(e^a_z(\mathbf{0},+\infty))=\mathbf{1}$. Indeed,
$e^a_z(\mathbf{0},+\infty)\geq e^a_z(\mathbf{0},t
\mathbf{1}]\succ\succ e^a_z(-\infty,\mathbf{0}]$. Observe that the preceding estimate immediately implies that
$$\mathbf{c}(e^a_z(\mathbf{0},+\infty))\geq
\mathbf{c}(e^a_z(-\infty,\mathbf{0}]).$$
By the definition of elements $e^a_z(-\infty,c]$, we have
$e^a_z(\mathbf{0},+\infty)\vee
e^a_z(-\infty,\mathbf{0}]=\mathbf{1}$. Consequently,
$$\mathbf{c}(e^a_z(\mathbf{0},+\infty))=\mathbf{c}(e^a_z(\mathbf{0},+\infty))\vee
\mathbf{c}(e^a_z(-\infty,\mathbf{0}])\geq
e^a_z(\mathbf{0},+\infty)\vee
e^a_z(-\infty,\mathbf{0}]=\mathbf{1}.$$

We want to show that if $d\in Z_+(LS(\mathcal{M}))$ and
$\mathbf{s}(d)=\mathbf{1}$ then $e^a_z(\mathbf{0},d]\succ\succ
e^a_z(-\infty,\mathbf{0}]+e^a_z(d,+\infty)$. Indeed,
a semiaxis $\mathbb{R}_+=(0,+\infty)$ can be split into countable
family of intervals $I_n=[\lambda_n,\mu_n)$, $n\geq1$. Then $e^d(I_n)\in
P(Z(\mathcal{M}))$, and since the spectral measure $e^d$ is $\sigma$-additive, we write

$$ \bigvee_{n=1}^\infty
e^d(I_n)=e^d(0,+\infty)=\mathbf{s}(d)=\mathbf{1}.$$
Next, assuming that $e^d(I_n)\neq 0$, we have

\begin{align*}
e^a_z(\mathbf{0},d]e^d(I_n) &=e^a_z(\mathbf{0},de^d(I_n)]e^d(I_n) \quad\text{(by Lemma \ref{l11} (ii) )}\\
&\geq e^a_z(\mathbf{0},\lambda_ne^d(I_n)]e^d(I_n) \quad\text{(by Lemma \ref{l11} (i) )}\\
&=e^a_z(\mathbf{0},\lambda_n\mathbf{1}]e^d(I_n) \quad\text{(by Lemma \ref{l11} (ii) )}\\
&\succ\succ
(e^a_z(-\infty,\mathbf{0}]+e^a_z(\lambda_n\mathbf{1},+\infty))e^d(I_n) \quad\text{(by the assumption)}\\
&=(e^a_z(-\infty,\mathbf{0}]+e^a_z(\lambda_ne^d(I_n),+\infty))e^d(I_n) \quad\text{(by Lemma \ref{l11} (ii) )}\\
&\geq (e^a_z(-\infty,\mathbf{0}]+e^a_z(d,+\infty))e^d(I_n)
\quad\text{(by Lemma \ref{l11} (i) )}
\end{align*}

This inequality implies, in particular, that for such choice of  $d$ the projection
$e^a_z(\mathbf{0},d]$ is properly infinite.

Next our goal is to construct a decreasing to zero sequence
$$\{d_n\}_{n=0}^\infty\subset Z_+(LS(\mathcal{M})),\ d_0\leq
\mathbf{1},\ d_{n+1}\leq d_n/2,\ \mathbf{s}(d_n)= \mathbf{1}$$ for
all $n\in \mathbb{N}$ and two sequences
$\{p_n\}_{n=0}^\infty,\ \{q_n\}_{n=0}^\infty$ of
pairwise disjoint projections in $\mathcal{M}$, which satisfy
following conditions:

(i). $p_nq_m=0,\ ap_n=p_na,\ aq_n=q_na,\ p_n\sim q_n$ for every
$n,m\geq 0$.

(ii). $p_n\leq e^a_z(d_n,+\infty),\ q_n\leq
e^a_z(-\infty,\varepsilon d_n]$ for every $n\geq 0$ и $q_0\geq
e^a_z(-\infty,\mathbf{0}]$.

(iii). $\bigvee_{n=0}^\infty p_n \vee \bigvee_{n=0}^\infty q_n =
\mathbf{1}$.

For any projection $p\in P(\mathcal{M})$ there exists a unique
central projection $z$ such that $pz$ is a finite projection and
$p(\mathbf{1}-z)$ is properly infinite projection and
$\mathbf{c}(p)\geq\mathbf{1}-z$ (if $p$ is finite projection then
$z=\mathbf{1}$, otherwise the assertion follows from \cite[Proposition 6.3.7]{KR2}).
In this case, we have
$\mathbf{c}(p(\mathbf{1}-z))=\mathbf{c}(p)(\mathbf{1}-z)=\mathbf{1}-z$
\cite[Proposition 5.5.3]{KR}. Let $z_0\in P(Z(\mathcal{M}))$ be
such a projection for $e^a_z(-\infty,\mathbf{0}]$ and let $z_n$ be
a such projection for $e^a_z(\mathbf{1}/n,+\infty)$,
$n\in\mathbb{N}$. Since the sequence $e^a_z(\mathbf{1}/n,+\infty)$
is non-decreasing, it follows that the sequence
$\{z_n\}_{n=1}^\infty$ is non-increasing. In addition we have
$\mathbf{1} = (\mathbf{1}-z_0)+\bigwedge_{n=1}^\infty z_0z_n
+z_0(\mathbf{1}-z_1) +\bigvee_{n=1}^\infty z_0(z_n-z_{n+1})$.
Hence, it is sufficient to prove the assertion for reduced
algebras
 $\mathcal{M}[\bigwedge_{n=1}^\infty
z_0z_n],\ \mathcal{M}[z_0(\mathbf{1}-z_1)],\
\mathcal{M}[z_0(z_n-z_{n+1})]$ and $\mathcal{M}(\mathbf{1}-z_0)$.
It is sufficient to consider three following cases:

(a). The projection $e^a_z(-\infty,\mathbf{0}]$ is finite and all projections
$e^a_z(\lambda\mathbf{1},+\infty)$ are the same for every $\lambda>0$.
The algebra
$\mathcal{M}[\bigwedge_{n=1}^\infty z_0z_n]$ satisfies this condition. Note that
in this case the projection $e^a_z(\mathbf{0},+\infty)$ is
a supremum of non-decreasing sequence of finite projections
 $e^a_z(\mathbf{1}/n,+\infty)$.

(b). The projection $e^a_z(-\infty,\mathbf{0}]$ is finite and there exists
 $\lambda>0$ such that the projection
$e_z^a(\lambda\mathbf{1},+\infty)$ is properly infinite and
$\mathbf{c}(e_z^a(\lambda\mathbf{1},+\infty))=\mathbf{1}$. Algebras
$\mathcal{M}[z_0(z_n-z_{n+1})]$
(in this case $\lambda=1/(n+1)$) and
$\mathcal{M}[z_0(\mathbf{1}-z_1)]$ (in this case
$\lambda=1$) satisfy this condition.

(c). The projection  $e^a_z(-\infty,\mathbf{0}]$ is properly infinite
and $\mathbf{c}(e^a_z(-\infty,\mathbf{0}])=\mathbf{1}$. The algebra
$\mathcal{M}(\mathbf{1}-z_0)$ satisfies this condition.

We would like to show that there exists an element $d_0\in Z_+(\mathcal{M}),\
d_0\leq \mathbf{1}$ such that $e_z^a(-\infty,\mathbf{0}]\prec\prec
e_z^a(d_0,+\infty)$ и $\mathbf{s}(d_0)=\mathbf{1}$.

Consider the case (a). We shall use Lemma \ref{l20}. To
this end, we set $b=\mathbf{1},\ \mu_n=1/n,\
q=e^a_z(-\infty,\mathbf{0}]$. In this case the projection
$q=e^a_z(-\infty,\mathbf{0}]$ is finite and projection
$e^a_z(\mathbf{0},+\infty)$ is properly infinite and
$\mathbf{c}(e^a_z(\mathbf{0},+\infty))=\mathbf{1}$. Hence,
the assumptions of Lemma \ref{l20} hold. Thus, there exists an element
$0<d_0\in Z(\mathcal{M})$  such that $d_0\leq \mathbf{1},\
\mathbf{s}(d_0)=\mathbf{1}$ and $e^a_z(-\infty,\mathbf{0}]\prec\prec
e^a_z(d_0,+\infty)$.

In the case (b) we set $d_0=\min(1,\lambda)\mathbf{1}$. Then
$e_z^a(-\infty,\mathbf{0}]$ is a finite projection and
$e_z^a(d_0,+\infty)$ is a properly infinite projection. Hence,
$e_z^a(-\infty,\mathbf{0}]\prec\prec e_z^a(d_0,+\infty)$ and
$\mathbf{s}(d_0)=\mathbf{1}$.

In the case (c) we will use Lemma \ref{l20} again. Set
$b=\mathbf{1},\ \mu_n=1/n,\ q=e^a_z(-\infty,\mathbf{0}]$. We have
that $q=e^a_z(-\infty,\mathbf{0}]$ is a properly
infinite projection, $\mathbf{c}(q)=\mathbf{1}$ and $q\prec\prec
e^a_z(\mathbf{0},+\infty)$. Hence, the assumptions of Lemma \ref{l20}
hold. So, there exists an element $0<d_0\in
Z(\mathcal{M})$ such that $d_0\leq \mathbf{1},\
\mathbf{s}(d_0)=\mathbf{1}$ and $e^a_z(-\infty,\mathbf{0}]\prec\prec
e^a_z(d_0,+\infty)$.

This completes the construction of the element $d_0$. Let us now show that
there exists a sequence $d_n\in Z_+(\mathcal{M})$ such that
$d_n\leq d_{n-1}/2$, $\mathbf{s}(d_n)=\mathbf{1}$ and
$e_z^a(d_n,+\infty)\succ\succ e_z^a(d_{n-1},+\infty)$ for every
$n\in \mathbb{N}$.

Suppose that elements $d_1,...,d_n$ have been already constructed.

We are going to use Lemma \ref{l20} again. For this we set
$b=d_n,\ \mu_m=1/(2m),\ q=e^a_z(d_n,+\infty)$. In the case (a)
 $q=e^a_z(d_n,+\infty)$ is finite projection. In the cases (b) and (c)
$e^a_z(d_n,+\infty)$ is a properly infinite projection,
$\mathbf{c}(q)=\mathbf{1}$ and $q\prec\prec
e^a_z(\mathbf{0},+\infty)$ (we have shown this fact at the beginning
of the proof). The assumptions of Lemma \ref{l20} hold. Hence,
there exists $0<d_{n+1}\in Z(\mathcal{M})$ such that
$d_{n+1}\leq d_n/2,\ \mathbf{s}(d_{n+1})=\mathbf{1}$ and
$e_z^a(d_{n+1},+\infty)\succ\succ e_z^a(d_n,+\infty)$.

Thus, we have constructed the sequence
$\{d_n\}\subset Z_+(LS(\mathcal{M}))$ such that $d_{n+1}\leq d_n/2$ and
$e^a_z(d_{n+1},+\infty)\succ\succ e^a_z(d_n,+\infty),\
\mathbf{s}(d_n)=\mathbf{1}$ for every $n\in\mathbb{N}$. In
addition we have $e^a_z(d_0,+\infty)\succ\succ e^a_z(-\infty,\mathbf{0}]$.

Set $p_0=e^a_z(d_0,+\infty)$. There exists a projection $r\in
P(\mathcal{M})$ such that $e^a_z(-\infty,\mathbf{0}]\sim r<p_0$.
By the assumption and using the argument at the beginning of the proof, we have
$e^a_z(\mathbf{0},\varepsilon d_0]\succ\succ p_0$, and so, by Lemma \ref{l1}, it follows
that there exists a projection
$q^1_0\in\mathcal{M}$ such that $p_0-r\sim q^1_0<
e^a_z(\mathbf{0},\varepsilon d_0]$ and $aq^1_0=q^1_0 a$ (in the case
(a) the condition (i) of Lemma \ref{l1} is applied and in the cases (b) and (c)
the condition (ii) is used). Set $q_0:=e^a_z(-\infty,\mathbf{0}]+q^1_0$.
Then $q_0\sim p_0$.

Suppose that projections $p_0,...,p_n;\ q_0,...,q_n$ have been constructed.
Set $p_{n+1}=e^a_z(d_{n+1},+\infty)
\prod_{k=0}^n(\mathbf{1}-p_k)\prod_{k=0}^n(\mathbf{1}-q_k)$. In the case
(a) all projections $p_k,\ q_k$ with $k\leq n$ are finite and
$e^a_z(\mathbf{0},\varepsilon d_{n+1}]$ is a properly infinite projection
which is a supremum of non-decreasing sequence of finite projections
$\{e^a_z(\mathbf{1}/m,\varepsilon d_{n+1}]\}_{m=1}^\infty$.
Hence, $e^a_z(\mathbf{0},\varepsilon d_n]
\prod_{k=0}^n(\mathbf{1}-p_k) \prod_{k=0}^n(\mathbf{1}-q_k)$ is
a properly infinite projection. It follows from Lemma \ref{l1} (i) that
there exists a projection $q_{n+1}\in\mathcal{M}$ such that
$p_{n+1}\sim q_{n+1}<e^a_z(\mathbf{0},\varepsilon d_{n+1}]
\prod_{k=0}^n(\mathbf{1}-p_k)\prod_{k=0}^n(\mathbf{1}-q_k)$ and
$aq_{n+1}=q_{n+1}a$.

Let now consider the cases (b) and (c). Recall that in these cases all $e^a_z(d_n,+\infty)$ are properly infinite projections.
Since $\sum_{k=0}^n p_k \leq
e^a_z(d_n,+\infty) \prec\prec e^a_z(d_{n+1},+\infty)$, by Lemma
\ref{l0} (ii) we obtain that $\sum_{k=0}^n p_k + \sum_{k=0}^n q_k
\prec\prec e^a_z(d_{n+1},+\infty)$. Then
$p_{n+1}=e^a_z(d_{n+1},+\infty)
\prod_{k=0}^n(\mathbf{1}-p_k)\prod_{k=0}^n(\mathbf{1}-q_k)=e^a_z(d_{n+1},+\infty)(\mathbf{1}-\sum_{k=0}^n p_k-\sum_{k=0}^n q_k)$ is a
properly infinite projection. It follows from Lemma \ref{l0} (iv) that
$p_{n+1}\sim e^a_z(d_{n+1},+\infty)\prec\prec
e^a_z(\mathbf{0},\varepsilon d_{n+1}]\sim
e^a_z(\mathbf{0},\varepsilon d_{n+1}]
\prod_{k=0}^n(\mathbf{1}-p_k)\prod_{k=0}^n(\mathbf{1}-q_k)$ (we applied Lemma \ref{l0} (iv) at the beginning and at the end of the chain).
So, it follows from Lemma \ref{l1} (ii) that there exists a projection
$q_{n+1}\in P(\mathcal{M})$ such that $q_{n+1} <
e^a_z(\mathbf{0},\varepsilon d_{n+1}]
\prod_{k=0}^n(\mathbf{1}-p_k)\prod_{k=0}^n(\mathbf{1}-q_k)$,
$q_{n+1}\sim p_{n+1}$ and $aq_{n+1}=q_{n+1}a$.

Thus, projections $p_{n+1}$ and $q_{n+1}$ are constructed.

It is clear that these projections satisfy conditions (i) and (ii). To check
the condition (iii) we note that
$\bigvee_{k=0}^n p_k\vee \bigvee_{k=0}^n q_k\geq
e^a_z(-\infty,\mathbf{0}]+e^a_z(d_n,+\infty) \uparrow \mathbf{1}$
with $n\rightarrow\infty$.

Now, we can proceed with the construction of the unitary operator
$u_\varepsilon\in\mathcal{M}$ from the assertion.

Let $v_n\in \mathcal{M}$ be a partial isometry such that
$v_n^*v_n=p_n,\ v_nv_n^*=q_n,\ n=0,1,...$. We set
\begin{align*} u_\varepsilon
=\sum_{n=0}^\infty v_n + \sum_{n=0}^\infty v_n^*
\end{align*} (here, the sums are taken in the strong operator topology).

Then, we have
\begin{align*} u_\varepsilon^*u_\varepsilon=\sum_{n=0}^\infty p_n +
\sum_{n=0}^\infty q_n = \mathbf{1},\ u_\varepsilon
u_\varepsilon^*=\sum_{n=0}^\infty q_n + \sum_{n=0}^\infty p_n =
\mathbf{1}.\end{align*}
Observe that
\begin{align*}u_\varepsilon p_n = q_n u_\varepsilon,\
u_\varepsilon q_n = p_n u_\varepsilon,\ ap_n=p_na,\ q_na=aq_n,\ n\ge 0,
\end{align*}
and so
the element $u_\varepsilon^* a u_\varepsilon$ commutes
with all the projections $p_n$ and $q_n$, $n\ge 0$. Moreover, since
for all $n\ge 0$, it holds
\begin{align*} ap_n=ae^a_z(d_n,+\infty)p_n\geq d_n
e^a_z(d_n,+\infty)p_n =d_n p_n,\\
aq_n=ae^a_z(-\infty,\varepsilon d_n)q_n\leq
\varepsilon d_n e^a_z(-\infty,\varepsilon d_n)q_n =
\varepsilon d_n q_n,
\end{align*}
we obtain immediately for all such $n$'s that
\begin{align*}
u_\varepsilon^* a
u_\varepsilon p_n &=u_\varepsilon^* aq_n u_\varepsilon \leq
\varepsilon d_n u_\varepsilon^* q_n
u_\varepsilon=\varepsilon d_n p_n,\\
u_\varepsilon^* a
u_\varepsilon q_n&=u_\varepsilon^* ap_n u_\varepsilon \geq
d_n u_\varepsilon^* p_n u_\varepsilon=d_n q_n.
\end{align*}
In particular, $(u_\varepsilon^* a u_\varepsilon - a)p_n\leq
\varepsilon d_n p_n-d_n p_n=
-d_n(1-\varepsilon)p_n\leq 0$. Taking into account that
$ap_n\geq d_n p_n$, we now obtain

\begin{align*}
|u_\varepsilon^* a
u_\varepsilon - a|p_n &=(a-u_\varepsilon^* a u_\varepsilon)p_n \geq
ap_n-\varepsilon d_n p_n\\ &\geq ap_n-\varepsilon
ap_n=(1-\varepsilon)ap_n\\ &=(1-\varepsilon)|a|p_n.
\end{align*}
Analogously, for every $n\ge 0$, we have $(u_\varepsilon^* a u_\varepsilon - a)q_n\geq d_n
q_n-\varepsilon d_n q_n= (1-\varepsilon) d_n q_n\geq 0$.
Therefore,
\begin{align*}
|u_\varepsilon^* a u_\varepsilon - a|q_n=(u_\varepsilon^*
a u_\varepsilon - a)q_n& \geq (1-\varepsilon) d_n q_n \\ &\geq
(1-\varepsilon)a q_n.
\end{align*}
Observe that the inequalities above hold for all $n\ge 0$. If
$n>0$, then $q_n<e^a_z(0,\varepsilon d_n]$, $q_na=aq_n$ by the
construction and so $aq_n=|a|q_n$, that is we have
$$
|u_\varepsilon^* a u_\varepsilon - a|q_n\geq (1-\varepsilon)|a|q_n.
$$

A little bit more care is required when $n=0$. In this case,
recall that $q_0=e^a_z(-\infty,\mathbf{0}]+q_0^1$, where
$q_0^1<e^a_z(\mathbf{0},\varepsilon d_0]$. Obviously,
$ae^a_z(-\infty,\mathbf{0}]\leq 0$, and so
$ae^a_z(-\infty,\mathbf{0}]=-|a|e^a_z(-\infty,\mathbf{0}]$. Therefore since (see above)
$u_\varepsilon^* a u_\varepsilon q_0\geq d_0 q_0$ and
$aq_0=ae^a_z(-\infty,\mathbf{0}]+aq_0^1=-|a|e^a_z(-\infty,\mathbf{0}]+aq_0^1$, we have
\begin{align*}
|u_\varepsilon^* a u_\varepsilon - a|q_0 & \geq (u_\varepsilon^* a
u_\varepsilon - a)q_0\geq  d_0
q_0-aq_0^1+|a|e^a_z(-\infty,\mathbf{0}]\\ & \geq
d_0q_0^1-\varepsilon d_0
q_0^1+|a|e^a_z(-\infty,\mathbf{0}]=(1-\varepsilon) d_0
q_0^1+|a|e^a_z(-\infty,\mathbf{0}]\\  & \geq (1-\varepsilon)a
q_0^1+|a|e^a_z(-\infty,\mathbf{0}]=
(1-\varepsilon)|a|q_0^1+|a|e^a_z(-\infty,\mathbf{0}]\\ &\geq
(1-\varepsilon)(|a|q_0^1+|a|e^a_z(-\infty,\mathbf{0}])=(1-\varepsilon)|a|q_0.
\end{align*}

Collecting all preceding inequalities, we see that for every $k\ge 0$ we have

\begin{align*}
|u_\varepsilon^* a u_\varepsilon - a|\sum_{n=0}^k (p_n+q_n)\geq (1-\varepsilon)|a|\sum_{n=0}^k (p_n+q_n)
\end{align*}
and since $\sum_{n=0}^\infty (p_n+q_n)=\mathbf{1}$, we conclude
\begin{align*}
|u_\varepsilon^* a u_\varepsilon - a|\geq (1-\varepsilon)|a|.
\end{align*}

The assertion of the lemma now follows by observing that
$|u_\varepsilon^* a u_\varepsilon -
a|=|[a,u_\varepsilon]|$.
\end{proof}

\textit{Proof of Theorem~\ref{main}.}
Prior to Lemma \ref{l33} we have shown that the identity $\mathbf{1}\in\mathcal{M}$ can be written as a sum
of three central projections
$\mathbf{1}=p_0+p_-+p_+$. The assertion (i) of Theorem~\ref{main} holds for the element $ap_0$ (by Lemma \ref{l3}) affiliated with the algebra
$\mathcal{M}p_0$. The assumptions of Lemma \ref{l2} hold for the element $(a-c_0)p_-$ in the algebra $\mathcal{M}p_-$.
Hence, the assertion (ii) of Theorem~\ref{main} holds in this algebra.
The assumptions of Lemma \ref{l2} hold for the element $(c_0-a)p_+$ in the algebra $\mathcal{M}p_+$ as well
(see the discussion preceding Lemma \ref{l1}). Hence, the assertions (ii) of Theorem~\ref{main} holds in this algebra as well.

Next, if $\mathcal{M}$ is finite or purely infinite
$\sigma$-finite algebra, then by Lemmas \ref{l31} and \ref{l32} we
have $\mathbf{1}\in P_0$, which implies $\mathbf{1}\leq p_0$. In
other words, we have $p_0=\mathbf{1}$ and by Lemma \ref{l3} the
assertions (i) of Theorem holds in this case.
\begin{flushright}$\Box$\end{flushright}

\end{document}